\documentclass[11pt]{amsart}
\usepackage[french]{babel}
\newtheorem{theoreme}{Th\'eor\`eme}[section]
\newtheorem{lemma}[theoreme]{Lemme}
\newtheorem{corollary}[theoreme]{Corollaire}
\newtheorem{definition}[theoreme]{D\'efinition}

\newcommand{\Oe}{\emptyset}
 \usepackage[latin1]{inputenc}
  \usepackage[T1]{fontenc}
  \usepackage{amsmath,amssymb}
  \usepackage[all]{xy}
\begin{document}
\title{Sur l'existence d'une cat\'egorie ayant une matrice donn\'ee}
\author{Samer Allouch}
\address{Laboratoire J. A. Dieudonn\'e\\Universit\'e de Nice-Sophia Antipolis}
\thanks{Ce papier a b\'en\'efici\'e d'une aide de l'Agence Nationale de la Recherche
portant la r\'ef\'erence ANR-09-BLAN-0151-02 (HODAG). }
\maketitle
\section{Introduction}
Ce papier est une continuation de \cite{Allouch2}, 
dont le but est d'\'etudier les correspondances entre
les cat\'egories finies d'ordre $n$ et les matrices carr\'ees de
taille $n$. Cette correspondance figure dans les papiers recents
de Leinster et Berger
\cite{Leinster} \cite{BergerLeinster}, voir aussi Kapranov \cite{Kapranov}.
L'\'etude des cat\'egories finies continue actuellement avec le papier
de Fiore, L\"uck et Sauer \cite{FioreLuckSauer}.

La question abord\'ee ici est de savoir, pour une matrice donn\'ee $M$,
s'il existe ou non une cat\'egorie $A$ associ\'ee \`a $M$. 
Dans \cite{Allouch2} il a \'et\'e trait\'e le cas des matrices strictement positives. 
Ici on consid\`ere une matrice $M$ quelconque. Apr\`es un rappel sur la condition
d'une matrice r\'eduite, la d\'efinition d'une relation $\geq$ sur les indices de la
matrice qui correspondent aux objets de la cat\'egorie
(la matrice est dite acceptable si la r\'elation est transitive et reflexive), et la partition de cet ensemble
en classes d'\'equivalence not\'ees $\lambda$,  les objets dans
chaque classe seront not\'es par $\lambda ^i$. Les classes $\lambda$ sont de type $U$ si $M(\lambda ^0,\lambda ^0)=1$ ou de type $V$ si $M(\lambda ^i,\lambda ^i)\geq 2$ pour tout $i$. 

Le r\'esultat du th\'eor\`eme \ref{gen1} est que
pour $M$ une matrice reduite alors
\begin{displaymath}
Cat(M)\neq{\emptyset}\Leftrightarrow \left\{
\begin{array}{llllll}
\qquad \qquad M  & \mbox{acceptable} &  \\
M(\lambda^{i},\lambda^{i})\geq & M(\lambda^{i},\lambda ^0) M(\lambda ^0,\lambda^{i}) +1& \forall \lambda \in U , i\geq 1\\
M(\lambda^{i},\lambda^{j})\geq & M(\lambda^{i},\lambda ^0) M(\lambda ^0,\lambda^{j}) & \forall \lambda \in U \\
M(\lambda^{i},\mu^{j}) \geq & M(\lambda^{i},\mu^{0}) & \forall \lambda > \mu , \mu \in U  \\
M(\lambda^{i},\mu^{j}) \geq & M(\lambda^{0},\mu^{j}) & \forall \lambda > \mu , \lambda \in U  \\
M(\lambda^{i},\mu^{j}) \geq  & M(\lambda^{0},\mu^{j})+
M(\lambda^{i},\mu^{0})
& \\
& \qquad -M(\lambda^{0},\mu^{0}) &\forall \lambda \geq \mu \in U
\end{array} \right.
\end{displaymath}
Pour expliquer le r\'esultat on traitera d'abord des exemples repr\'esentatifs,
y compris le cas des matrices strictement positives de \cite{Allouch2} pour fixer
les notations.

\section{Partition des matrices}
Soit M=$(a_{ij})$ une matrice carr\'ee d'order n, $A$ une cat\'egorie associ\'ee \`a  $M$ dont les objets sont  $Ob(A)=\{x_1,....,x_n\}$ et $A(x_i,x_j)= a_{ij}$ avec $A(x_i,x_j)$ l'ensemble des morphismes de $$f:x_i\longmapsto x_j$$.\\
On d\'efinit $\sim$ une relation de semi-\'equivalence sur $Ob(A)$
par:
\begin{displaymath}
x_i \sim x_j\Longleftrightarrow\left\{ \begin{array}{l}
a_{ij}>0\\
a_{ji}>0
\end{array} \right.
\end{displaymath}
On va partitionner $Ob(A)$ par $A/_\sim$=$\{\alpha_1,...,\alpha_q\}$=$\{U_1,.....U_p,V_1,.....,V_q\}$ o$\grave{u}$ $U_i$ contient au moins $x_j$ telque $a_{jj}=1$  pour tout $i\in\lbrace 1,...,p\rbrace$\\
et $V_j$ ne contient aucune $x_s$ telque $a_{ss}=1$.\\
Soient $\lambda, \mu, \varphi \ldots \in Ob(A)/\sim $ sont des
classes objets. Les objets d'une classe $\lambda$ seront not\'es
$\lambda ^i$ pour $1\leq i\leq |\lambda |$ si $\lambda \in V$, ou
$0\leq i < |\lambda |$ si $\lambda \in U$. De cette fa\c{c}on, les
objets qui n'ont que l'identit\'e comme endomorphismes, sont les
$\lambda ^0$.\\
Ensuite, les morphismes seront not\'es $\lambda ^i \mu ^j X^k$
o\`u $X$ d\'esigne une lettre majuscule, eventuellement avec un
exposant $k$ qui pourrait eventuellement \^etre une paire.

Le choix de lettre d\'esignera le type de morphisme.

On commence avec les notations pour les morphismes $\lambda
^i\lambda ^jX$.

D'abord, pour l'identit\'e on \'ecrira $\lambda ^i\lambda ^iI$
(pas besoin d'exposant car il n'y a qu'une seule identit\'e pour
chaque objet $\lambda ^i$).

Ensuite, on aura des morphismes de la forme $\lambda ^i\lambda ^j
F^{u,v}$ avec les conditions suivantes: on a
$$
1\leq u \leq a(\lambda ^i), \;\; 1\leq v \leq b(\lambda ^j)
$$
avec (si $\lambda \in U$) $a(\lambda ^0)=1$ et $b(\lambda ^0)= 1$;
en particulier pour $i=0$ il n'y a que $u=1$, pour $j=0$ il n'y a
que $v=1$. On fera la convention que
$$
\lambda ^0\lambda ^0F^{1,1} = \lambda ^0\lambda ^0I
$$
est l'identit\'e; cependant pour $i>0$ on a
$$
\lambda ^i\lambda ^iF^{1,1} \neq \lambda ^i\lambda ^iI.
$$

Pour $\lambda \in U$ et soit $i=0$, soit $j=0$, ces morphismes
sont les seuls morphismes.

Pour $i,j\geq 1$, et dans tous les cas $\lambda \in V$, on peut
avoir en plus des morphismes de la forme $\lambda ^i\lambda ^j
G^k$, pour
$$
1\leq k \leq M(\lambda ^i ,\lambda ^j) - a(\lambda ^i) b(\lambda
^j) \mbox{   si  } i\neq j
$$
ou pour
$$
1\leq k \leq M(\lambda ^i ,\lambda ^j) - a(\lambda ^i) b(\lambda
^j) -1\mbox{   si  } i= j.
$$
Ce nombre de morphismes suppl\'ementaires peut etre \'egale \`a $0$,
dans ce cas il n'y en a pas. Ces $G^k$ s'occupent
de la technique ``d'ajouter des
morphismes''.

 On consid\`ere maintenant les morphismes dans le cas $\lambda
^i\mu ^j$ avec $\lambda >\mu$, et en supposant par exemple que
$\lambda ,\mu \in U$. Il y aura des types de morphismes suivants:
$\lambda ^i\mu ^j A^k$, $\lambda ^i\mu ^j B^k$, $\lambda ^i\mu ^j
C^k$, $\lambda ^i\mu ^j D^k$

Pour  $\lambda ^i\mu ^j A^k$,
$$
1\leq k \leq M(\lambda ^0,\mu ^0).
$$

Pour $\lambda ^i\mu ^j B^k$,
$$
1\leq k \leq M(\lambda ^i,\mu ^0) - M(\lambda ^0,\mu ^0).
$$

Pour $\lambda ^i\mu ^j C^k$,
$$
1\leq k \leq M(\lambda ^0,\mu ^j) - M(\lambda ^0,\mu ^0).
$$

Pour $\lambda ^i\mu ^j D^k$,
$$
1\leq k \leq M(\lambda ^i,\mu ^j) - M(\lambda ^i,\mu ^0) -
M(\lambda ^0,\mu ^j) + M(\lambda ^0,\mu ^0).
$$
\textbf{Exemple}: Soit M une matrice d\'efinie par:
\begin{displaymath} \mathbf{M} = \left(
\begin{array}{cc|cc}
1& b&c&d  \\
e&f&k&l\\
\hline
0&0&1&x\\
0&0&q&m
\end{array} \right) .
\end{displaymath}
avec a,b,c,d,e,f,k,l,x,q,m sont strictement positives.\\
alors si cat(M)$\neq {\Oe}$ donc $A/_\sim$=$\{\lambda,\mu\}$ avec
$\lambda=\{x_1,x_2\}$ et $\mu=\{x_3,x_4\}$.
\section {Matrice reduite}
\begin{definition}:
Soit $A$ une categorie d'ordre $n$ avec objets $x_1,\ldots , x_n$,
on dit que $x_i$ et $x_j$ sont {\em isomorphes} s'il existe $f\in
A(x_i,x_j)$ et $g\in A(x_j,x_i)$ tels que $fg = 1_{x_j}$ et $gf =
1_{x_i}$.
\end{definition}
 \textbf{Rq}: Si $x_i$ et $x_j$ sont isomorphes, alors pour tout objet $x_k$ on a des isomorphismes
d'ensembles
$$
A(x_k,x_i)\stackrel{\cong}{\rightarrow} A(x_k,x_j),
$$
donn\'es par $h\mapsto fh$ dans une direction, et $u\mapsto gu$
dans l'autre; et
$$
A(x_i,x_k)\stackrel{\cong}{\rightarrow} A(x_j,x_k),
$$
donn\'e par $h\mapsto hg$ dans une direction, et $u\mapsto uf$
dans l'autre. Si $M$ est la matrice de $A$, on en d\'eduit:
$$
\forall k, \;\; M_{ki} = M_{kj}$$et$$\forall k, \;\; M_{ik} =
M_{jk}.
$$ \\
\begin{definition}:Soit
\label{defreduite} $A$ est une cat\'egorie telle qu'il existe deux
objets distincts $x_i$ et $x_j$ ($i\neq j$) qui sont isomorphes,
on dira que $A$ est {\em non-r\'eduite}. On dira que $A$ est {\em
r\'eduite} sinon,c'est-\`a-dire si deux objets distincts sont
toujours non-isomorphes. On dira qu'une matrice $M$ est {\em
non-r\'eduite} s'il existe $i\neq j$ tel que$$\forall k, \;\;
M_{ki} = M_{kj}$$et$$\forall k, \;\; M_{ik} = M_{jk},$$cela veut
dire que la ligne $i$ \'egale la ligne $j$ et la colonne $i$
\'egale la colonne $j$. On dira qu'une matrice $M$ est {\em
r\'eduite} si elle n'est pas non-r\'eduite.
\end{definition}
\textbf{Rq:} D'apres le debut ci-dessus, on obtient que si $A$ est
non-r\'eduite, alors $M$ est non-r\'eduite. Donc, par contrapos\'e
si $M$ est r\'eduite alors $A$ est r\'eduite. Le contraire n'est
pas forc\'emment vrai: il peut exister une cat\'egorie $A$ telle
que $M$ est non-r\'eduite, mais $A$ r\'eduite, par exemple on peut
avoir une cat\'egorie $A$ d'ordre $2$ dont la matrice
non-r\'eduite est
$$
M=\left(
\begin{array}{cc}2 & 2 \\2 & 2\end{array}\right)
$$
mais telle que les deux objets de $A$ sont non-isomorphes et donc
$A$ r\'eduite.
\begin{theoreme}: Si $M$ une matrice non r\'eduite, on peut r\'eduire
$M$ en  une sous matrice $N$ r\'eduite telle que M marche si et
seulement si  $N$ marche.
\end{theoreme}
En effet: Supposons que $M$ est une matrice $n\times n$
non-r\'eduite. On peut d\'efinir une r\'elation d'\'equivalence
sur l'ensemble d'indices $\{ 1,\ldots ,n\}$ en disant que $i\sim
j$ si $\forall k, \;\; M_{ki} = M_{kj}$ et$\forall k, \;\; M_{ik}
= M_{jk}$. Celle-ci est sym\'etrique, reflexive et transitive. On
obtient donc une partition de l'ensemble d'indices en r\'eunion
disjointe de sous-ensembles$$\{ 1,\ldots ,n\} = U_1\sqcup U_2
\sqcup \cdots \sqcup U_m$$avec $U_a\cap U_b= \emptyset$, telle que
tous les \'el\'ements d'un $U_a$ donn\'e sont \'equivalents, et
les \'el\'ements de $U_a$ ne sont pas \'equivalents aux
\'el\'ements de $U_b$ pour $a\neq b$. (Ce sont les classes
d'\'equivalence pour la r\'elation d'\'equivalence). Choisissons
un repr\'esantant $r(a)\in U_a$ pour chaque classe
d'\'equivalence. Dans l'autre sens, on note par $c(i)\in \{
1,\ldots , m\}$ l'unique \'el\'ement telle que $i\in U_{c(i)}$.
Ici $c(i)$ est la classe d'\'equivalence contenant $i$. On
a$$c(r(a))=a$$mais $r(c(i))$ n'est pas toujours \'egale \`a $i$:
on a seulement qu'ils sont \'equivalents $r(c(i))\sim i$.  On
obtient une sous-matrice de taille $m\times m$
$$
N_{ab}:= M_{r(a),r(b)}.
$$
Il est a noter qu'on peut faire en sorte que $r(a)<r(b)$ pour
$a<b$: on choisit $r(a)$ le plus petit \'el\'ement de $U_a$, et on
num\'erote les classes $U_a$ par ordre croissant de leur plus
petit \'el\'ement. Dans ce cas $N$ est vraiementune sous-matrice
de $M$. Noter que $N$ est r\'eduite, puisque les \'el\'ements de
$U_a$ et $U_b$ ne sont pas \'equivalents pour $a\neq b$. Si $A$
est une cat\'egorie dont la matrice est $M$, on obtient une
sous-cat\'egorie pleine $B\subset A$ qui consiste des objets
$r(a)$ seulement, $a=1,\ldots , m$. La matrice de $B$ est $N$. On
conclut que si $M$ marche, alors $N$ marche. L'\'equivalence entre
$i$ et $r(c(i))$ implique que pour tout $k$ on a$$M_{k,i} =
M_{k,r(c(i))},\;\;\; M_{i,k} = M_{r(c(i)),k}.$$ On en d\'eduit que
pour tout $i,j$ on a$$M_{i,j}= M_{r(c(i)),j} = M_{r(c(i)),r(c(j))}
= N_{c(i),c(j)}.$$Ceci indique comment aller dans l'autre
sens.Supposons que $B$ est une cat\'egorie dont la matrice est
$N$. Notons par $y_1,\ldots , y_m$ les objets de $B$. On d\'efinit
une cat\'egorie $A$ avec objets not\'es $x_1,\ldots , x_n$ en
posant$$A(x_i,x_j) \cong B(y_{c(i)},y_{c(j)}).$$Si on veut etre
plus precise on pourrait d\'efinir$$A(x_i,x_j) := \{ (i,j,\beta),
\;\;\; \beta \in B(y_{c(i)},y_{c(j)}) \} .$$ La composition est la
m\^eme que celle de $B$, i.e.$$(i,j,\beta )(j,k,\beta '):=
(i,k,\beta \beta ').$$De m\^eme pour les identit\'es, et
l'associativit\'e et les r\`egles des identit\'es sont faciles \`a
v\'erifier. Donc $A$ est une cat\'egorie.\\
On a:
$$
|A(x_i,x_j)| = |B(y_{c(i)},y_{c(j)})| = N_{c(i),c(j)} = M_{i,j}.
$$
Donc $A$ corr\'espond \`a la matrice $M$. De cette discussion on
conclut: \'etant donn\'ee une matrice non-r\'eduite $M$, on peut
construire par la construction pr\'ec\'edente une sous-matrice $N$
qui est r\'eduite, telle que $M$ marche si et seulement si $N$
marche. La sous-matrice $N$ est unique \`a permutation d'indices
pr\`es.
\begin{lemma}:Soit $M$ est une matrice r\'eduite avec $M_{i,j}>0$, et s'il existe
$i\neq j$ tels que $M_{i,i}=1$ et $M_{j,j}=1$, alors $M$ ne marche
pas.
\end{lemma}
En effet: On suppose que $M$ marche alors il existe une
cat\'egorie $A$ associ\'ee \`a $M$ et comme $M$ est r\'eduite
alors $A$ est r\'eduite. En plus $M_{i,i}=1$ et $M_{j,j}=1$ alors
$x_i $ et $x_j$ sont isomorphes.  En effet, $A(x_i,x_j)$ a
$M_{ij}>0$ \'el\'ements, on peut en choisir un $f$; et
$A(x_j,x_i)$ a $M_{ji}>0$ \'el\'ements, choisissons-en $g$. Alors
$fg=1$ et $gf=1$ car $|A(x_i,x_i)|=M_{ii}=1$ et
$|A(x_j,x_j)|=M_{jj}=1$. Alors $A$ est non-r\'eduite contradiction
donc $M$
ne marche pas.\\
\section {Matrices stictement positives}

Les sections 4 \`a 6, qui constituent un rappel de \cite{Allouch2}, 
expliqueront la m\'ethode ``d'ajouter des morphismes'' qui sera
ensuite utilis\'ee sans plus d'explication pour la preuve du th\'eor\`eme
\ref{gen1}. 

Soit M=($a_{ij}$) une matrice carre\'e d$'$ordre n a coefficients
positives.\\
\begin{theoreme}[Leinster \cite{BergerLeinster}]
\label{leinsterthm}: Soit $M=(m_{ij})$ une matrice carr\'ee dont
les coificients sont des entiers naturels et pour tout i
$m_{ii}\geq 2$ , alors Cat(M)$\neq{\Oe}$ (i.e.d il existe une
cat\'egorie associ\'e \`a M).
\end{theoreme}
En effet: Soit $M=(m_{ij})$ de taille avec $m_{ii}\geq 2$, on pose
$n_{ij}:= m_{ij}$ pour $i\neq j$ et $n_{ii}:= m_{ii}-1$. On peut
d\'efinir une semi-cat\'egorie A associ\'e \`a $N$ dont les objets
sont $1,2,.....,n$, pour tout couple $(i,j)$ on a une fl\`eche
$\Phi_{ij}$ :i $\rightarrow$ j tel que $\Phi_{ij} \neq 1_{ii}$,la
loi de composition d\'efinit par si $f:i \rightarrow j$ et
$g:j\rightarrow k$ $\Phi_{ij}$ alors $gf=\Phi_{ik}$. Ensuite on
peut d\'efinir une cat\'egorie $B$ en rajoutant \`a $A$ les
identit\'es, pour tout $i$ on a $1_{ii}: i \rightarrow i$.  La
matrice de $B$ est $M$.\\
\\
Soit M une matrice carr\'ee d'order 2 d\'efinie
par:\begin{displaymath} \mathbf{M} = \left( \begin{array}{cc}
1 & b  \\
c & d
\end{array} \right)
\end{displaymath}
avec $b,c,d>1$.
\begin{theoreme}
Cat(M)$\neq{\Oe}$$ \Leftrightarrow$ d$\geq bc+1$
\end{theoreme}

Pour $d=bc+1$\\
$ \Rightarrow\big)$:on a Cat(M)$\neq{\Oe}$ On va démontrer que  d$\geq bc+1$\\
Soit A une cat\'egorie associ\'ee \`a M ,Comme tous les
coefficientes de la matrice sont strictement positives alors il y
a une seule classe \'equivante not\'e
$\lambda$.\\
Soit A une cat\'egorie associe à M dont les objets sont
$\lambda$=$\{\lambda^0,\lambda^1\}$.\\
les morphismes sont d\'efinis par:\\
\\
$A(\lambda^0\lambda^0)$=${I}$.\\
$A(\lambda^1\lambda^1)$=$\{\lambda^1\lambda^1G^{u,v}\}$\\
$A(\lambda^0\lambda^1)$=$\{\lambda^0\lambda^1G^{u,v} /1\leq u\leq
a(\lambda^0)=1$, $1\leq v \leq
b(\lambda^1)=b\}$=$\{\lambda^0\lambda^1G^{1,v}/1\leq v \leq b\}$\\
$A(\lambda^1\lambda^0)$=$\{\lambda^1\lambda^0G^{u,v} /1\leq u\leq
a(\lambda^1)=c$, $1\leq v \leq
b(\lambda^0)=1\}$=$\{\lambda^1\lambda^0G^{u,1}/1\leq u \leq c\}$\\
Ona:$(\lambda^1\lambda^0G^{u,1})(\lambda^0\lambda^1G^{1,v})$=$I_0$
pour tout
u,v\\\\
\textbf{Rq (1)}:il n'existe pas u,v telque
$(\lambda^0\lambda^1G^{1,v})(\lambda^1\lambda^0G^{u,1})$=$I_1$,
sinon s'il existe u,v telque
$(\lambda^0\lambda^1G^{1,v})(\lambda^1\lambda^0G^{u,1})$=$I_1$,
on a
$(\lambda^1\lambda^0G^{u,1})\big[(\lambda^1\lambda^1G^{n,m})(\lambda^0\lambda^1G^{1,v})\big]$=I$_0$
alors
$(\lambda^0\lambda^1G^{1,v})\big[(\lambda^1\lambda^0G^{u,1})(\lambda^1\lambda^1G^{n,m})(\lambda^0\lambda^1G^{1,v})\big]$=$(\lambda^0\lambda^1G^{1,v})$\\
donc
$\big[(\lambda^0\lambda^1G^{1,v})\big(\lambda^1\lambda^0G^{u,1})\big]\big[(\lambda^1\lambda^1G^{n,m})(\lambda^0\lambda^1G^{1,v})\big]$=$(\lambda^0\lambda^1G^{1,v})$\\
alors $(\lambda^1\lambda^1G^{n,m})(\lambda^0\lambda^1G^{1,v})$=$(\lambda^0\lambda^1G^{1,v})$\\
d'autre
part,$\big[(\lambda^1\lambda^1G^{n,m})(\lambda^0\lambda^1G^{1,v})\big]$($(\lambda^1\lambda^0G^{u,1})$=$(\lambda^0\lambda^1G^{1,v})(\lambda^1\lambda^0G^{u,1})$=I$_1=(\lambda^1\lambda^1G^{n,m})$
contarduction.\\
\\
\textbf{Rq (2)}:soient
$(\lambda^0\lambda^1G^{1,v})(\lambda^1\lambda^0G^{u,1})$=$(\lambda^1\lambda^1G^{n,m})$
alors
$(\lambda^1\lambda^1G^{n,m})^2$=$(\lambda^1\lambda^1G^{n,m}).$\\
En effet
$(\lambda^1\lambda^0G^{u,1})(\lambda^0\lambda^1G^{1,v})$=$I_0$\\
alors
$(\lambda^0\lambda^1G^{1,v})\big[(\lambda^1\lambda^0G^{u,1})(\lambda^0\lambda^1G^{1,v})\big]$=$(\lambda^0\lambda^1G^{1,v})$\\
donc $\big[(\lambda^0\lambda^1G^{1,v})(\lambda^1\lambda^0G^{u,1})\big](\lambda^0\lambda^1G^{1,v})$=$(\lambda^0\lambda^1G^{1,v})$\\
alors $(\lambda^1\lambda^1G^{n,m})(\lambda^0\lambda^1G^{1,v})$=$(\lambda^0\lambda^1G^{1,v})$\\
d'autre part
$\big[(\lambda^1\lambda^1G^{n,m})(\lambda^0\lambda^1G^{1,v})\big](\lambda^1\lambda^0G^{u,1})$=$(\lambda^0\lambda^1G^{1,v})(\lambda^1\lambda^0G^{u,1})$
=$(\lambda^1\lambda^1G^{n,m})$\\
$(\lambda^1\lambda^1G^{n,m})\big[(\lambda^0\lambda^1G^{1,v})(\lambda^1\lambda^0G^{u,1})\big]$=$(\lambda^1\lambda^1G^{n,m})^2$\\\\
\textbf{Rq (3)}:s'il existe
$(\lambda^0\lambda^1G^{1,v})(\lambda^1\lambda^0G^{u,1})$=$(\lambda^0\lambda^1G^{1,p})(\lambda^1\lambda^0G^{u,1})$
alors p=v.\\
En effet
$\big[(\lambda^0\lambda^1G^{1,v})(\lambda^1\lambda^0G^{u,1})\big](\lambda^0\lambda^1G^{1,p})$\\
=$(\lambda^0\lambda^1G^{1,v})\big[(\lambda^1\lambda^0G^{u,1})(\lambda^0\lambda^1G^{1,p})\big]$\\
=$(\lambda^0\lambda^1G^{1,v})$\\
=$\big[(\lambda^0\lambda^1G^{1,p})(\lambda^1\lambda^0G^{u,1})\big](\lambda^0\lambda^1G^{1,p})$\\
=$(\lambda^0\lambda^1G^{1,p})\big[(\lambda^1\lambda^0G^{u,1})(\lambda^0\lambda^1G^{1,p})\big]$\\
=$(\lambda^0\lambda^1G^{1,p})$\\
Donc $(\lambda^0\lambda^1G^{1,v})$=$(\lambda^0\lambda^1G^{1,p})$.\\\\
\textbf{Rq (3)}:s'il existe
$(\lambda^0\lambda^1G^{1,v})(\lambda^1\lambda^0G^{u,1})$=$(\lambda^0\lambda^1G^{1,v})(\lambda^1\lambda^0G^{p,1})$
alors p=u.\\
En effet
$\big[(\lambda^1\lambda^0G^{p,1})(\lambda^0\lambda^1G^{1,v})\big](\lambda^1\lambda^0G^{u,1})$\\
=$(\lambda^1\lambda^0G^{u,1})$\\
=$(\lambda^1\lambda^0G^{p,1})\big[(\lambda^0\lambda^1G^{1,v})(\lambda^1\lambda^0G^{u,1})\big]$\\
=$(\lambda^1\lambda^0G^{p,1})\big[(\lambda^0\lambda^1G^{1,v})(\lambda^1\lambda^0G^{p,1})\big]$\\
=$\big[(\lambda^1\lambda^0G^{p,1})(\lambda^0\lambda^1G^{1,v})\big](\lambda^1\lambda^0G^{p,1})$\\
=$(\lambda^1\lambda^0G^{p,1})$\\
donc =$(\lambda^1\lambda^0G^{u,1})$==$(\lambda^1\lambda^0G^{p,1})$\\\\
\textbf{Rq (3)}:il n'existe pas $u\neq m$ et $v\neq n$ telque
$(\lambda^0\lambda^1G^{1,v})(\lambda^1\lambda^0G^{u,1})$=$(\lambda^0\lambda^1G^{1,n})(\lambda^1\lambda^0G^{m,1})$.\\
En effet
$\big[(\lambda^1\lambda^0G^{m,1})(\lambda^0\lambda^1G^{1,v})\big](\lambda^1\lambda^0G^{u,1})$\\
=$(\lambda^1\lambda^0G^{u,1})$\\
=$(\lambda^1\lambda^0G^{m,1})\big[(\lambda^0\lambda^1G^{1,v})(\lambda^1\lambda^0G^{u,1})\big]$\\
=$(\lambda^1\lambda^0G^{m,1})\big[(\lambda^0\lambda^1G^{1,n})(\lambda^1\lambda^0G^{m,1})\big]$\\
=$(\lambda^1\lambda^0G^{m,1})$\\
donc m=u contradiction .\\
Donc les 4 remarques \textbf{Rq (1)} \textbf{Rq (3)} \textbf{Rq (4)} et \textbf{Rq (5)} donnent $(d-bc)>0$ c.\`a.d $d>bc$\\
\\
$\Leftarrow\big)d>bc\Rightarrow Cat(M)\neq{\Oe}.$ \\
Pour $d=bc+1$\\
soit A une cat\'egorie dont les objets sont $\lambda^0,\lambda^1$
et les morphismes sont d\'efinis comme
pr\'ec\`edant .\\
On va d\'efinir la loi de composition par:\\
$(\lambda^j\lambda^pG^{u,v})(\lambda^i\lambda^jG^{n,m})$=$(\lambda^i\lambda^pG^{n,v})$
pour tout i,j,p$\in\{0,1\}$c.\`a.d:\\
$(\lambda^1\lambda^0G^{u,1})(\lambda^0\lambda^1G^{1,v})=(\lambda^0\lambda^0G^{1,1}).$\\
$(\lambda^1\lambda^1G^{u,v})(\lambda^0\lambda^1G^{1,v'})=(\lambda^0\lambda^1G^{1,v}).$\\
$(\lambda^0\lambda^1G^{1,v})(\lambda^1\lambda^0G^{u,1})=(\lambda^1\lambda^1G^{u,v}).$\\
$(\lambda^1\lambda^0G^{u',1})(\lambda^1\lambda^1G^{u,v})=(\lambda^1\lambda^0G^{u,1}).$\\
$(\lambda^1\lambda^1G^{u,v})(\lambda^1\lambda^1G^{u',v'})=(\lambda^1\lambda^1G^{u',v}).$\\
soient i,j,p,q $\in\{0,1\}$ alors:\\
$\big[(\lambda^p\lambda^qG^{x,y})(\lambda^j\lambda^pG^{u,v})\big](\lambda^i\lambda^jG^{n,m})$=\\
$(\lambda^j\lambda^qG^{u,y})(\lambda^i\lambda^jG^{n,m})$=\\
$(\lambda^i\lambda^qG^{n,y})$\\
d'autre part\\
$(\lambda^p\lambda^qG^{x,y})\big[(\lambda^j\lambda^pG^{u,v})(\lambda^i\lambda^jG^{n,m})\big]$=\\
$(\lambda^p\lambda^qG^{x,y})(\lambda^i\lambda^pG^{n,v})=$\\
$(\lambda^i\lambda^qG^{n,y})$\\
Donc
$\big[(\lambda^p\lambda^qG^{x,y})(\lambda^j\lambda^pG^{u,v})\big](\lambda^i\lambda^jG^{n,m})$=$(\lambda^p\lambda^qG^{x,y})\big[(\lambda^j\lambda^pG^{u,v})(\lambda^i\lambda^jG^{n,m})\big]$
ce qui donne A cat\'egorie associ\'ee \`a M donc $Cat(M)\neq{\Oe}.$ \\
Pour $d>$bc+1 alors d=bc+1+n\\
Soit A$'$ une cat\'egorie tel que Ob(A)= Ob(A$'$)tel que: \\
$A'$($\lambda^1\lambda^1$)=$\{\lambda^1\lambda^1G^{u,v} / 1\leq u\leq b $ et $1\leq v\leq c\}\cup\{\lambda^1\lambda^1K^p/1\leq p\leq n\}$ \\
$A'(\lambda^0\lambda^0)$=${I}$.\\
$A'(\lambda^0\lambda^1)$=$\{\lambda^0\lambda^1G^{u,v} /1\leq u\leq
a(\lambda^0)=1$, $1\leq v \leq
b(\lambda^1)=b\}$=$\{\lambda^0\lambda^1G^{1,v}/1\leq v \leq b\}$\\
$A'(\lambda^1\lambda^0)$=$\{\lambda^1\lambda^0G^{u,v} /1\leq u\leq
a(\lambda^1)=c$, $1\leq v \leq
b(\lambda^0)=1\}$=$\{\lambda^1\lambda^0G^{u,1}/1\leq u \leq c\}$\\
avec la loi de composition d\'efinie par:\\
$(\lambda^1\lambda^1K^i)(\lambda^0\lambda^1G^{1,v})=(\lambda^1\lambda^1G^{1,1})(\lambda^0\lambda^1G^{1,v})=(\lambda^0\lambda^1G^{1,1})$\\
$(\lambda^1\lambda^0G^{u,1})(\lambda^1\lambda^1K^i)=(\lambda^1\lambda^0G^{u,1})(\lambda^1\lambda^1G^{b,c})=(\lambda^1\lambda^0G^{b,1})$\\
$(\lambda^1\lambda^1G^{u,v})(\lambda^1\lambda^1K^i)=(\lambda^1\lambda^1G^{u,v})(\lambda^1\lambda^1G^{b,c})=(\lambda^1\lambda^1G^{b,v})$\\
$(\lambda^1\lambda^1K^i)(\lambda^1\lambda^1G^{u,v})=(\lambda^1\lambda^1G^{1,1})(\lambda^1\lambda^1G^{u,v})=(\lambda^1\lambda^1G^{u,1})$\\
$(\lambda^1\lambda^1K^i)(\lambda^1\lambda^1K^{i'})=(\lambda^1\lambda^1G^{1,1})(\lambda^1\lambda^1G^{b,c})=(\lambda^1\lambda^1G^{b,1})$\\

\begin{corollary}
soit $M$ une matrice carr\'e d'order 2 d\'efinie par:
\begin{displaymath}
\mathbf{M} = \left( \begin{array}{cc}
a & b  \\
c & d
\end{array} \right)
\end{displaymath}
avec $a,b,c,d$ sont strictement positives, alors $cat(M)\neq{\Oe}$ dans les cas suivantes :\\
1- a=b=c=d=1 voir \cite{Allouch}\\
2- a=1, $d>bc$.\\
3- d=1, $a>bc$.\\
4- $a>1$ , $d>1$ voir \cite{Leinster}.\\
si non $Cat(M)={\Oe}$
\end{corollary}
\section {Matrices triples}
Soit M une matrice d'order 3 strictement positives d\'efinie par:
\begin{displaymath}
\mathbf{M} = \left( \begin{array}{ccc}
1 & a &b  \\
c & n &m\\
p & q &r
\end{array} \right)
\end{displaymath}
$n>1$ et $r>1$, sinon la matrice n'est pas r\'eduite. On peut
d\'eduire de l'\'etude pr\'ec\`edant que $n\geq ac+1$ et $r\geq
bp+1$ en effet les deux matrices suivantes
\begin{displaymath}
\left( \begin{array}{cc}
1 & b  \\
p & r
\end{array} \right)
\end{displaymath}
et
\begin{displaymath}
\left( \begin{array}{cc}
1 & a  \\
c & n
\end{array} \right)
\end{displaymath}
sont sous-matices voir \cite {Allouch} et si N est une sous matrice qui ne
marche pas alors M ne marche pas.

En plus $m\geq bc$ et $q\geq ap$ en effet:\\
soit $A$ une cat\'egorie associ\'ee \`a M dont les objets sont $\{x_1,x_2,x_3\}$\\
D'apr\'es la decomposition de matrice, comme M est strictement
positives alors il y a une seule bloc c.\`a.d il y a une seule
classe $\lambda$ et les objets d\'efinis par
$\{\lambda^0,\lambda^1,\lambda^2\}$
 et les morphismes sont donn\'es par:\\
$A(\lambda^0\lambda^0)$=$\{id_{\lambda^0}=1\}$\\
$A(\lambda^2\lambda^1)$=$\{\lambda^2\lambda^1G^{u,v}\}$\\
$A(\lambda^1\lambda^2)$=$\{\lambda^1\lambda^0G^{u,v}\}$\\
$A(\lambda^0\lambda^1)$=$\{\lambda^0\lambda^1G^{1,v}/1\leq v \leq a\}$\\
$A(\lambda^1\lambda^0)$=$\{\lambda^1\lambda^0G^{u,1}/1\leq u \leq c\}$\\
$A(\lambda^0\lambda^2)$=$\{\lambda^0\lambda^2G^{1,v}/1\leq v \leq b\}$\\
$A(\lambda^2\lambda^0)$=$\{\lambda^2\lambda^0G^{u,1}/1\leq u \leq p\}$\\
$A(\lambda^1\lambda^1)$=$\{\lambda^1\lambda^1G^{u,v}/1\leq u \leq a, 1\leq v \leq c\}$\\
$A(\lambda^2\lambda^2)$=$\{\lambda^2\lambda^2G^{u,v}/1\leq u \leq b, 1\leq v \leq p\}$\\
Nous revenons au but: il faut d\'emontrer que $m\geq bc$ et $q\geq ap$ \\
on suppose que $m<bc$ alors il y a 3 cas: \\
a)il existe  u$\neq u'$,v$\neq v'$,x avec  $1\leq v,v'\leq b$ et
$1\leq u,u'\leq c$ et $1\leq x\leq p$ tel que :\\
$(\lambda^0\lambda^2G^{1,v})(\lambda^1\lambda^0G^{u,1})=(\lambda^0\lambda^2G^{1,v'})(\lambda^1\lambda^0G^{u',1})$\\
alors\\
$(\lambda^2\lambda^0G^{x,1})\big[(\lambda^0\lambda^2G^{1,v})(\lambda^1\lambda^0G^{u,1})\big]$=
$(\lambda^2\lambda^0G^{x,1})\big[(\lambda^0\lambda^2G^{1,v'})(\lambda^1\lambda^0G^{u',1})\big]$\\
donc\\
$\big[(\lambda^2\lambda^0G^{x,1})(\lambda^0\lambda^2G^{1,v})\big](\lambda^1\lambda^0G^{u,1})$=
$\big[(\lambda^2\lambda^0G^{x,1})(\lambda^0\lambda^2G^{1,v'})\big](\lambda^1\lambda^0G^{u',1})$\\
alors\\
$id_{\lambda^0}(\lambda^1\lambda^0G^{u,1})$=$id_{\lambda^0}(\lambda^1\lambda^0G^{u',1})$\\
alors\\
$u=u'$ contradiction donc n'existe pas cette cas .\\
\\
b) il existe u,v$\neq v'$,x avec  $1\leq u,u'\leq c$ $1\leq v,v'\leq b$ et $1\leq x\leq a$ tel que :\\
$(\lambda^0\lambda^2G^{1,v})(\lambda^1\lambda^0G^{u,1})=(\lambda^0\lambda^2G^{1,v'})(\lambda^1\lambda^0G^{u,1})$\\
alors\\
$\big[(\lambda^0\lambda^2G^{1,v})(\lambda^1\lambda^0G^{u,1})\big](\lambda^0\lambda^1G^{1,x})$=
$\big[(\lambda^0\lambda^2G^{1,v'})(\lambda^1\lambda^0G^{u,1})\big](\lambda^0\lambda^1G^{1,x})$\\
donc\\
$(\lambda^0\lambda^2G^{1,v})\big[(\lambda^1\lambda^0G^{u,1})(\lambda^0\lambda^1G^{1,x})\big]$=
$(\lambda^0\lambda^2G^{1,v'})\big[(\lambda^1\lambda^0G^{u,1})(\lambda^0\lambda^1G^{1,x})\big]$\\
alors\\
$(\lambda^0\lambda^2G^{1,v})id_{\lambda^0}$=$(\lambda^0\lambda^2G^{1,v'})id_{\lambda^0}$\\
alors\\
$v=v'$ contradiction donc n'existe pas cette cas .\\
\\
c)il existe  u$\neq u'$,v,x avec  $1\leq v\leq b$ et
$1\leq u,u'\leq c$ et $1\leq x\leq p$ tel que :\\
$(\lambda^0\lambda^2G^{1,v})(\lambda^1\lambda^0G^{u,1})=(\lambda^0\lambda^2G^{1,v})(\lambda^1\lambda^0G^{u',1})$\\
alors\\
$(\lambda^2\lambda^0G^{x,1})\big[(\lambda^0\lambda^2G^{1,v})(\lambda^1\lambda^0G^{u,1})\big]$=
$(\lambda^2\lambda^0G^{x,1})\big[(\lambda^0\lambda^2G^{1,v})(\lambda^1\lambda^0G^{u',1})\big]$\\
donc\\
$\big[(\lambda^2\lambda^0G^{x,1})(\lambda^0\lambda^2G^{1,v})\big](\lambda^1\lambda^0G^{u,1})$=
$\big[(\lambda^2\lambda^0G^{x,1})(\lambda^0\lambda^2G^{1,v})\big](\lambda^1\lambda^0G^{u',1})$\\
alors\\
$id_{\lambda^0}(\lambda^1\lambda^0G^{u,1})$=$id_{\lambda^0}(\lambda^1\lambda^0G^{u',1})$\\
alors\\
$u=u'$ contradiction donc n'existe pas cette cas .\\
Donc $m\geq bc$.\\
De la meme pour $q\geq ap$\\
\begin{theoreme}
Soit M une matrice triple dont les coeficientes sont strictement
positive d\'efinie par:
\begin{displaymath}
\mathbf{M} = \left( \begin{array}{ccc}
1 & a &b  \\
c & n &m\\
p & q &r
\end{array} \right)
\end{displaymath}
alors $n=ac+1,r=bp+1,m=bc,q=ap \Rightarrow Cat(M)\neq \emptyset$
\end{theoreme}
En effet:soit A une semi-cat\'egorie d'order 3 dont les objets
sont $\{\lambda^0\lambda^1,\lambda^2\}$ et les morphismes d\'efinis par:\\
$A(\lambda^0\lambda^0)$=$\{id_{\lambda^0}=1\}$\\
$A(\lambda^0\lambda^1)$=$\{\lambda^0\lambda^1G^{1,v}/1\leq v \leq a\}$\\
$A(\lambda^1\lambda^0)$=$\{\lambda^1\lambda^0G^{u,1}/1\leq u \leq c\}$\\
$A(\lambda^0\lambda^2)$=$\{\lambda^0\lambda^2G^{1,v}/1\leq v \leq b\}$\\
$A(\lambda^2\lambda^0)$=$\{\lambda^2\lambda^0G^{u,1}/1\leq u \leq p\}$\\
$A(\lambda^2\lambda^1)$=$\{\lambda^2\lambda^1G^{u,v}/1\leq u \leq p,1\leq v \leq a\}$\\
$A(\lambda^1\lambda^2)$=$\{\lambda^1\lambda^0G^{u,v}/1\leq u \leq c,1\leq v \leq b\}$\\
$A(\lambda^1\lambda^1)$=$\{\lambda^1\lambda^1G^{u,v}\neq 1/1\leq u \leq a, 1\leq v \leq c\}$\\
$A(\lambda^2\lambda^2)$=$\{\lambda^2\lambda^2G^{u,v}\neq 1/1\leq u \leq b, 1\leq v \leq p\}$\\
\`a partir de ces \'equations $n=ac+1,r=bp+1,m=bc$ et $q=ap $ nous
pouvons d\'efinir la loi de
composition par:\\
$(\lambda^j\lambda^kG^{u',v'})(\lambda^i\lambda^jG^{u,v})$=$(\lambda^i\lambda^kG^{u,v'})\forall
i,j,k \in \{0,1,2\}$\\
on va v\'erifier l'associativit\'e:\\
$\Big[(\lambda^k\lambda^sG^{u'',v''})(\lambda^j\lambda^kG^{u',v'})\Big](\lambda^i\lambda^jG^{u,v})=$\\
$(\lambda^j\lambda^sG^{u',v''})(\lambda^i\lambda^jG^{u,v})=$\\
$(\lambda^i,\lambda^sG^{u,,v''})$\\
d$'autre$ part
$(\lambda^k\lambda^sG^{u'',v''})\Big[(\lambda^j\lambda^kG^{u',v'})(\lambda^i\lambda^jG^{u,v})\Big]=$\\
$(\lambda^k\lambda^sG^{u'',v''})(\lambda^i\lambda^kG^{u,,v'})=$\\
$(\lambda^i,\lambda^sG^{u,,v''})$\\
donc
$\Big[(\lambda^k\lambda^sG^{u'',v''})(\lambda^j\lambda^kG^{u',v'})\Big](\lambda^i\lambda^jG^{u,v})=$
$(\lambda^k\lambda^sG^{u'',v''})\Big[(\lambda^j\lambda^kG^{u',v'})(\lambda^i\lambda^jG^{u,v})\Big]$\\
alors A est une semi-cat\'egorie associ\'e \`a M.\\
Soit $B=A\oplus \{id_{x_2}\}\oplus\{id_{x_3}\}$ donc $B$ est une cat\'egorie associ\'e \`a $M$ ce qui donne $Cat(M)\neq {\Oe} $ \\
\textbf{Notation:}\\
On a pour la matrice triple $M$ d\'efinie par:
\begin{displaymath}
\mathbf{M} = \left( \begin{array}{ccc}
1 & a &b  \\
c & n &m\\
p & q &r
\end{array} \right)
\end{displaymath}
avec $n=ac$ ,$r=bp$, $m=bc$ et $q=ap$ d'apres pr\'ec\`edant $M$ admet $A$ comme semi-cat\'egorie. \\
\\
Maintenant on va chercher une semi-cat\'egorie associ\'e \`a $M$ avec $n>ac,r>bp,m>bc , q>ap $ et apr\`es on ajoute les identit\'es. \\
Soit $M(ac+1)$ matrice d\'efinit par :\\
\begin{displaymath}
\mathbf{M(ac+1)} = \left( \begin{array}{ccc}
1 & a    &b  \\
c & ac+1 &m\\
p & q    &r
\end{array} \right)
\end{displaymath}
avec $r=bp,m=bc ,q=ap $\\
soit $A'$ une semi-cat\'egorie dont les objets sont $Ob(A')=Ob(A)$  avec les morphismes sont: \\
$A'(\lambda^0\lambda^0)$=$\{id_{\lambda^0}=1\}$\\
$A'(\lambda^0\lambda^1)$=$\{\lambda^0\lambda^1G^{1,v}/1\leq v \leq a\}$\\
$A'(\lambda^1\lambda^0)$=$\{\lambda^1\lambda^0G^{u,1}/1\leq u \leq c\}$\\
$A'(\lambda^0\lambda^2)$=$\{\lambda^0\lambda^2G^{1,v}/1\leq v \leq b\}$\\
$A'(\lambda^2\lambda^0)$=$\{\lambda^2\lambda^0G^{u,1}/1\leq u \leq p\}$\\
$A'(\lambda^2\lambda^1)$=$\{\lambda^2\lambda^1G^{u,v}/1\leq u \leq p,1\leq v \leq a\}$\\
$A'(\lambda^1\lambda^2)$=$\{\lambda^1\lambda^0G^{u,v}/1\leq u \leq c,1\leq v \leq b\}$\\
$A'(\lambda^2\lambda^2)$=$\{\lambda^2\lambda^2G^{u,v}\neq 1/1\leq u \leq b, 1\leq v \leq p\}$\\
$A'(\lambda^1\lambda^1)$=$\{\lambda^1\lambda^1G^{u,v}\neq 1/1\leq u \leq a, 1\leq v \leq c\}\bigcup \{\lambda^1\lambda^1K^1\}$\\
Les \'equations de la loi composition de $A'$ sont les memes que de $A$ en plus les \'equations dependant de $e'$ sont:\\
$(\lambda^1\lambda^1K^1)(\lambda^i\lambda^1G^{u,v})$=$(\lambda^1\lambda^1G^{1,1})(\lambda^i\lambda^1G^{u,v})$=$(\lambda^i\lambda^1G^{u,1})$\\
$(\lambda^1\lambda^iG^{u,v})(\lambda^1\lambda^1K^1)$=$(\lambda^1\lambda^iG^{u,v})(\lambda^1\lambda^1G^{a,c})$=$(\lambda^1\lambda^iG^{a,v})$\\
$(\lambda^1\lambda^1K^1)(\lambda^1\lambda^1K^1)$=$(\lambda^1\lambda^1K^1)$\\
Pour l$'$associativit\'e il y a six possibilit\'es des \'equations
d\'efinies par:\\
$(\lambda^1\lambda^jG^{u,v})\Big[(\lambda^1\lambda^1K^1)(\lambda^i\lambda^1G^{u',v'})\Big]$=\\
$(\lambda^1\lambda^jG^{u,v})\Big[(\lambda^1\lambda^1G^{1,1})(\lambda^i\lambda^1G^{u',v'})\Big]$=\\
$(\lambda^1\lambda^jG^{u,v})(\lambda^i\lambda^1G^{u',1})$=\\
$(\lambda^i\lambda^jG^{u',v})$\\
D$'$autre part \\
$\Big[(\lambda^1\lambda^jG^{u,v})(\lambda^1\lambda^1K^1)\Big](\lambda^i\lambda^1G^{u',v'})$=\\
$\Big[(\lambda^1\lambda^jG^{u,v})(\lambda^1\lambda^1G^{a,c})\Big](\lambda^i\lambda^1G^{u',v'})$=\\
$(\lambda^1\lambda^jG^{a,v})(\lambda^i\lambda^1G^{u',v'})$=\\
$(\lambda^i\lambda^jG^{u',v})$\\
Donc $(\lambda^1\lambda^jG^{u,v})\Big[(\lambda^1\lambda^1K^1)(\lambda^i\lambda^1G^{u',v'})\Big]$=$\Big[(\lambda^1\lambda^jG^{u,v})(\lambda^1\lambda^1K^1)\Big](\lambda^i\lambda^1G^{u',v'})$\\
$(\lambda^1\lambda^1K^1)\Big[(\lambda^j\lambda^1G^{u,v})(\lambda^i\lambda^jG^{u',v'})\Big]$=\\
$(\lambda^1\lambda^1K^1)(\lambda^i\lambda^1G^{u',v})$=\\
$(\lambda^1\lambda^1G^{1,1})(\lambda^i\lambda^1G^{u',v})$=\\
$(\lambda^i\lambda^1G^{u',1})$\\
$\Big[(\lambda^1\lambda^1K^1)(\lambda^j\lambda^1G^{u,v})\Big](\lambda^i\lambda^jG^{u',v'})$=\\
$\Big[(\lambda^1\lambda^1K^1)(\lambda^j\lambda^1G^{u,v})\Big](\lambda^i\lambda^jG^{u',v'})$=\\
$(\lambda^1\lambda^1G^{1,1})(\lambda^j\lambda^1G^{u,v})\Big](\lambda^i\lambda^jG^{u',v'})$=\\
$(\lambda^j\lambda^1G^{u,1})(\lambda^i\lambda^jG^{u',v'})$=\\
$(\lambda^i\lambda^1G^{u',,1})$\\
Alors $(\lambda^1\lambda^1K^1)\Big[(\lambda^j\lambda^1G^{u,v})(\lambda^i\lambda^jG^{u',v'})\Big]$=$\Big[(\lambda^1\lambda^1K^1)(\lambda^j\lambda^1G^{u,v})\Big](\lambda^i\lambda^jG^{u',v'})$\\
$(\lambda^i\lambda^iG^{u,v})\Big[(\lambda^1\lambda^iG^{u',v'})(\lambda^1\lambda^1K^1)\Big]$=\\
$(\lambda^i\lambda^iG^{u,v})\Big[(\lambda^1\lambda^iG^{u',v'})(\lambda^1\lambda^1G^{a,c})\Big]$=\\
$(\lambda^i\lambda^iG^{u,v})(\lambda^1\lambda^iG^{a,v'})$=\\
$(\lambda^1\lambda^iG^{a,v})$
$\Big[(\lambda^i\lambda^iG^{u,v})(\lambda^1\lambda^iG^{u',v'})\Big](\lambda^1\lambda^1K^1)$=\\
$\Big[(\lambda^i\lambda^iG^{u,v})(\lambda^1\lambda^iG^{u',v'})\Big](\lambda^1\lambda^1G^{a,c})=$\\
$(\lambda^1\lambda^iG^{u',v})(\lambda^1\lambda^1G^{a,c})=$\\
$(\lambda^1\lambda^iG^{a,v})$ \\
Donc $(\lambda^i\lambda^iG^{u,v})\Big[(\lambda^1\lambda^iG^{u',v'})(\lambda^1\lambda^1K^1)\Big]$=$\Big[(\lambda^i\lambda^iG^{u,v})(\lambda^1\lambda^iG^{u',v'})\Big](\lambda^1\lambda^1K^1)$\\
$\Big[(\lambda^1\lambda^1K^1)(\lambda^1\lambda^1K^1)\Big](\lambda^1\lambda^1K^1)$=$(\lambda^1\lambda^1K^1)\Big[(\lambda^1\lambda^1K^1)(\lambda^1\lambda^1K^1)\Big]=(\lambda^1\lambda^1K^1)$\\
$(\lambda^1\lambda^1K^1)^2(\lambda^i\lambda^1G^{u,v})=$\\
$(\lambda^1\lambda^1K^1)(\lambda^i\lambda^1G^{u,v})=$\\
$(\lambda^i\lambda^1G^{u,1})=$\\
$(\lambda^1\lambda^1K^1)\Big[(\lambda^1\lambda^1K^1)(\lambda^i\lambda^1G^{u,v})\Big]=$\\
$(\lambda^1\lambda^1K^1)\Big[(\lambda^1\lambda^1G^{1,1})(\lambda^i\lambda^1G^{u,v})\Big]=$\\
$(\lambda^1\lambda^1K^1)(\lambda^i\lambda^1G^{u,1})$\\
$(\lambda^i\lambda^1G^{u,1})=$\\
Alors $(\lambda^1\lambda^1K^1)^2(\lambda^i\lambda^1G^{u,v})=(\lambda^1\lambda^1K^1)\Big[(\lambda^1\lambda^1K^1)(\lambda^i\lambda^1G^{u,v})\Big]$\\
$(\lambda^i\lambda^1G^{u,v})\Big[(\lambda^1\lambda^1K^1)(\lambda^1\lambda^1K^1)\Big]$=\\
$(\lambda^i\lambda^1G^{u,v})(\lambda^1\lambda^1K^1)$=\\
$(\lambda^1\lambda^1G^{a,v})$\\
$\Big[(\lambda^i\lambda^1G^{u,v})(\lambda^1\lambda^1K^1)\Big](\lambda^1\lambda^1K^1)$=\\
$(\lambda^1\lambda^1G^{a,v})(\lambda^1\lambda^1K^1)$=\\
$(\lambda^1\lambda^1G^{a,v})$\\
Donc $(\lambda^i\lambda^1G^{u,v})\Big[(\lambda^1\lambda^1K^1)(\lambda^1\lambda^1K^1)\Big]$=$\Big[(\lambda^i\lambda^1G^{u,v})(\lambda^1\lambda^1K^1)\Big](\lambda^1\lambda^1K^1)$\\
Donc $A_{1}$ est une semi-cat\'egorie associ\'e \`a M.\\
on suppose que $A_{n-1}$ semi-cat\'egorie tel que
$A_{n-1}=A'\cup\{(\lambda^1\lambda^1K^2),(\lambda^1\lambda^1K^3),...,(\lambda^1\lambda^1K^{n-1})\}
=A\cup\{(\lambda^1\lambda^1K^1),(\lambda^1\lambda^1K^2),...,(\lambda^1\lambda^1K^{n-1})\}$ avec $K^i\neq K^j$  $\forall i,j\in \{1,...,(n-1)\}$ \\
la loi composition d\'efinie  par:\\
$(\lambda^1\lambda^1K^i)(...)=(\lambda^1\lambda^1G^{1,1})(...)$\\
$(...)(\lambda^1\lambda^1K^i)=(...)(\lambda^1\lambda^1G^{a,c})$\\
$(\lambda^1\lambda^1K^i)(\lambda^1\lambda^1K^j)=(\lambda^1\lambda^1G^{a,1})\forall i\neq j$\\
$(\lambda^1\lambda^1K^i)(\lambda^1\lambda^1K^i)=(\lambda^1\lambda^1K^i)\forall i$\\
On va v\'erifier les associativit\'es, comme nous
allons \'etudier les \'equations par rapport \`a
$(\lambda^1\lambda^1K^1)$ alors il reste 4 \'equations pour
v\'erifier:\\
$\Big[(\lambda^1\lambda^1K^i)(\lambda^1\lambda^1K^j)\Big](\lambda^i\lambda^1G^{u,v})=$\\
$(\lambda^1\lambda^1G^{a,1})(\lambda^i\lambda^1G^{u,v})=$\\
$(\lambda^i\lambda^1G^{u,1})$\\
$(\lambda^1\lambda^1K^i)\Big[(\lambda^1\lambda^1K^j)(\lambda^i\lambda^1G^{u,v})\Big]=$\\
$(\lambda^1\lambda^1K^i)\Big[(\lambda^1\lambda^1G^{1,1})(\lambda^i\lambda^1G^{u,v})\Big]=$\\
$(\lambda^1\lambda^1K^i)(\lambda^i\lambda^1G^{u,1})=$\\
$(\lambda^i\lambda^1G^{u,1})$\\
Donc $\Big[(\lambda^1\lambda^1K^i)(\lambda^1\lambda^1K^j)\Big](\lambda^i\lambda^1G^{u,v})=$$(\lambda^1\lambda^1K^i)\Big[(\lambda^1\lambda^1K^j)(\lambda^i\lambda^1G^{u,v})\Big]$\\
$\Big[(\lambda^1\lambda^iG^{u,v})(\lambda^1\lambda^1K^i)\Big](\lambda^1\lambda^1K^j)=$\\
$\Big[(\lambda^1\lambda^iG^{u,v})(\lambda^1\lambda^1G^{a,c})\Big](\lambda^1\lambda^1K^j)=$\\
$(\lambda^1\lambda^iG^{a,v})(\lambda^1\lambda^1K^j)=$\\
$(\lambda^1\lambda^iG^{a,v})$\\
$(\lambda^1\lambda^iG^{u,v})\Big[(\lambda^1\lambda^1K^i)(\lambda^1\lambda^1K^j)\Big]=$\\
$(\lambda^1\lambda^iG^{u,v})(\lambda^1\lambda^1G^{a,1})=$\\
$(\lambda^1\lambda^iG^{a,v})$\\
Alors $\Big[(\lambda^1\lambda^iG^{u,v})(\lambda^1\lambda^1K^i)\Big](\lambda^1\lambda^1K^j)=(\lambda^1\lambda^iG^{u,v})\Big[(\lambda^1\lambda^1K^i)(\lambda^1\lambda^1K^j)\Big]$\\
$\Big[(\lambda^1\lambda^1K^i)(\lambda^1\lambda^1G^{u,v})\Big](\lambda^1\lambda^1K^j)=$\\
$\Big[(\lambda^1\lambda^1G^{1,1})(\lambda^1\lambda^1G^{u,v})\Big](\lambda^1\lambda^1K^j)=$\\
$(\lambda^1\lambda^1G^{u,1})(\lambda^1\lambda^1K^j)=$\\
$(\lambda^1\lambda^1G^{a,1})$\\
$(\lambda^1\lambda^1K^i)\Big[(\lambda^1\lambda^1G^{u,v})(\lambda^1\lambda^1K^j)\Big]=$\\
$(\lambda^1\lambda^1K^i)(\lambda^1\lambda^1G^{a,v})=$\\
$(\lambda^1\lambda^1G^{a,1})$\\
Donc $\Big[(\lambda^1\lambda^1K^i)(\lambda^1\lambda^1G^{u,v})\Big](\lambda^1\lambda^1K^j)=(\lambda^1\lambda^1K^i)\Big[(\lambda^1\lambda^1G^{u,v})(\lambda^1\lambda^1K^j)\Big]$\\
$\Big[(\lambda^1\lambda^1K^i)(\lambda^1\lambda^1K^j)\Big](\lambda^1\lambda^1K^p)$=\\
$(\lambda^1\lambda^1G^{a,1})(\lambda^1\lambda^1K^p)$=\\
$(\lambda^1\lambda^1G^{a,1})$\\
$(\lambda^1\lambda^1K^i)\Big[(\lambda^1\lambda^1K^j)(\lambda^1\lambda^1K^p)\Big]$=\\
$(\lambda^1\lambda^1K^i)(\lambda^1\lambda^1G^{a,1})$=\\
$(\lambda^1\lambda^1G^{a,1})$\\
Alors $\Big[(\lambda^1\lambda^1K^i)(\lambda^1\lambda^1K^j)\Big](\lambda^1\lambda^1K^p)$=$(\lambda^1\lambda^1K^i)\Big[(\lambda^1\lambda^1K^j)(\lambda^1\lambda^1K^p)\Big]$\\
Donc $A_{n-1}$ est une semi-cat\'egorie associe\'e \`a la matrice suivant:\\
\begin{displaymath}
\mathbf{M(ac+(n-1))} = \left( \begin{array}{ccc}
1 & a        &b  \\
c & ac+(n-1) &bc\\
p & ap       &bp
\end{array} \right)
\end{displaymath}
Maintenant on ajoute sur  $A(\lambda^1,\lambda^2)$ \\
soient $\{(\lambda^1\lambda^2H^1),(\lambda^1\lambda^2H^2),...,(\lambda^1\lambda^2H^{m})\}$ morphismes dans $A(x_2,x_3)$  avec $k^i\neq k^j$ pour tout $i,j\in \{1,...,m\}$ tel que le loi de composition d\'efinie par:\\
$(\lambda^1\lambda^2H^i)(...)=(\lambda^1\lambda^2G^{1,1})(...)$\\
$(...)(\lambda^1\lambda^2H^i)=(...)(\lambda^1\lambda^2G^{b,c})$\\
$(\lambda^1\lambda^2H^i)(\lambda^1\lambda^1K^j)=(\lambda^1\lambda^2G^{1,1})(\lambda^1\lambda^2G^{a,c})=(\lambda^1\lambda^2G^{a,1})$\\
les \'equations associatives associ\'ees \`a $(\lambda^1\lambda^2H^i)$ sont:\\
$(\lambda^1\lambda^2H^i)\Big[(\lambda^i\lambda^1G^{u,v})(\lambda^j\lambda^iG^{u',v'})\Big]$=\\
$(\lambda^1\lambda^2H^i)(\lambda^j\lambda^1G^{u',v})$=\\
$(\lambda^j\lambda^2G^{u',1})$\\
$\Big[(\lambda^1\lambda^2H^i)(\lambda^i\lambda^1G^{u,v})\Big](\lambda^j\lambda^iG^{u',v'})$=\\
$(\lambda^i\lambda^2G^{u,1})(\lambda^j\lambda^iG^{u',v'})$=\\
$(\lambda^j\lambda^2G^{u',1})$\\
Alors $\Big[(\lambda^1\lambda^2H^i)(\lambda^i\lambda^1G^{u,v})\Big](\lambda^j\lambda^iG^{u',v'})$=$(\lambda^1\lambda^2H^i)\Big[(\lambda^i\lambda^1G^{u,v})(\lambda^j\lambda^iG^{u',v'})\Big]$\\
$(\lambda^2\lambda^jG^{u,v})\Big[(\lambda^1\lambda^2H^i)(\lambda^i\lambda^1G^{u',v'})\Big]$=\\
$(\lambda^2\lambda^jG^{u,v})(\lambda^i\lambda^2G^{u',1})$=\\
$(\lambda^i\lambda^jG^{u',v})$\\
$\Big[(\lambda^2\lambda^jG^{u,v})(\lambda^1\lambda^2H^i)\Big](\lambda^i\lambda^1G^{u',v'})$=\\
$(\lambda^1\lambda^jG^{b,v})(\lambda^i\lambda^1G^{u',v'})$=\\
$(\lambda^i\lambda^jG^{u',v})$\\
Donc $(\lambda^2\lambda^jG^{u,v})\Big[(\lambda^1\lambda^2H^i)(\lambda^i\lambda^1G^{u',v'})\Big]$=$\Big[(\lambda^2\lambda^jG^{u,v})(\lambda^1\lambda^2H^i)\Big](\lambda^i\lambda^1G^{u',v'})$\\
$(\lambda^1\lambda^2H^i)\Big[(\lambda^i\lambda^1G^{u,v})(\lambda^j\lambda^iG^{u',v'})\Big]=$
$(\lambda^1\lambda^2H^i)(\lambda^j\lambda^1G^{u',v})=$\\
$(\lambda^j\lambda^2G^{u',1})$\\
$\Big[(\lambda^1\lambda^2H^i)(\lambda^i\lambda^1G^{u,v})\Big](\lambda^j\lambda^iG^{u',v'})=$\\
$(\lambda^i\lambda^2G^{u,1})(\lambda^j\lambda^iG^{u',v'})=$\\
$(\lambda^j\lambda^2G^{u',1})$\\
Alors
$(\lambda^1\lambda^2H^i)\Big[(\lambda^i\lambda^1G^{u,v})(\lambda^j\lambda^iG^{u',v'})\Big]=$$\Big[(\lambda^1\lambda^2H^i)(\lambda^i\lambda^1G^{u,v})\Big](\lambda^j\lambda^iG^{u',v'})$\\
$(\lambda^1\lambda^2H^i)\Big[(\lambda^1\lambda^1K^j)(\lambda^1\lambda^1K^p)\Big]$=\\
$(\lambda^1\lambda^2H^i)(\lambda^1\lambda^1G^{a,1}$=\\
$(\lambda^1\lambda^2G^{a,1}$=\\
$\Big[(\lambda^1\lambda^2H^i)(\lambda^1\lambda^1K^j)\Big](\lambda^1\lambda^1K^p)$=\\
$(\lambda^1\lambda^2G^{a,1})(\lambda^1\lambda^1K^p)$=\\
$(\lambda^1\lambda^2G^{a,1}$=\\
Donc $(\lambda^1\lambda^2H^i)\Big[(\lambda^1\lambda^1K^j)(\lambda^1\lambda^1K^p)\Big]$=$\Big[(\lambda^1\lambda^2H^i)(\lambda^1\lambda^1K^j)\Big](\lambda^1\lambda^1K^p)$\\
$(\lambda^1\lambda^2H^i)\Big[(\lambda^1\lambda^1K^j)(\lambda^i\lambda^1G^{u,v})\Big]$=\\
$(\lambda^1\lambda^2H^i)(\lambda^i\lambda^1G^{u,1})$=\\
$(\lambda^i\lambda^2G^{u,1})$=\\
$\Big[(\lambda^1\lambda^2H^i)(\lambda^1\lambda^1K^j)\Big](\lambda^i\lambda^1G^{u,v})$=\\
$(\lambda^1\lambda^2G^{a,1})(\lambda^i\lambda^1G^{u,v})$=\\
$(\lambda^i\lambda^2G^{u,1})$=\\
Alors $(\lambda^1\lambda^2H^i)\Big[(\lambda^1\lambda^1K^j)(\lambda^i\lambda^1G^{u,v})\Big]$=$\Big[(\lambda^1\lambda^2H^i)(\lambda^1\lambda^1K^j)\Big](\lambda^i\lambda^1G^{u,v})$\\
De la meme construction on peux ajouter aussi sur $A(x_3,x_3)$ et sur $A(x_3,x_2)$ des morphismes adjoints avec la d\'efinition de la loi de composition:\\
$(\lambda^2\lambda^2N^i)(...)=(\lambda^2\lambda^2G^{1,1})(...)$\\
$(...)(\lambda^2\lambda^2N^i)=(...)(\lambda^2\lambda^2G^{b,p})$\\
$(\lambda^2\lambda^2N^i)(\lambda^2\lambda^2N^j)=(\lambda^2\lambda^2G{b,1})$\\
$(\lambda^2\lambda^2N^i)^2$=$(\lambda^2\lambda^2N^i)$\\
$(\lambda^2\lambda^1M^i)(...)=(\lambda^2\lambda^1G^{1,1})(...)$\\
$(...)(\lambda^2\lambda^1M^i)=(...)(\lambda^2\lambda^1G^{a,p})$\\
$(\lambda^1\lambda^1K^i)(\lambda^2\lambda^1M^j)=(\lambda^2\lambda^1G^{a,1})$\\
$(\lambda^2\lambda^1M^i)(\lambda^2\lambda^2N^j)=(\lambda^2\lambda^1G^{1,p})$\\
$(\lambda^1\lambda^2H^i)(\lambda^2\lambda^1M^j)=(\lambda^2\lambda^2G^{a,1})$\\
$(\lambda^2\lambda^1M^i)(\lambda^1\lambda^2M^j)=(\lambda^1\lambda^1M^{b,1})$\\
Les \'equations associatives marchent , donc $A_{(n-1,m,q,r-1)}$ est une semi-cat\'egorie  \\
c.\`a .d  en ajoutant les identites $A'_{(n,m,q,r)}$ est une cat\'egorie  associ\'ee \`a la matrice $M$ qui est d\'efinie par:\\
\begin{displaymath}
\mathbf{M} = \left( \begin{array}{ccc}
1 & a        &b\\
c & ac+n     &bc+m\\
p & ap+q       &bp+r
\end{array} \right)
\end{displaymath}
o\`u $n,r,m,q$ sont des entiers naturels.\\
\begin{corollary}
Soit $M$ une matrice d'order 3 tel que :
\begin{displaymath}
\mathbf{M} = \left( \begin{array}{ccc}
z & a        &b\\
c & n     &m\\
p & q       &r
\end{array} \right)
\end{displaymath}\\
avec $z\geq 1$,$n>ac$ ,$r>bp$, $m\geq bc$ et $q \geq ap$ alors $Cat(M)\neq {\Oe}$\\
\end{corollary}
\textbf{Preuve}: soit $N$ une matrice d\'efinie par:
\begin{displaymath}
\mathbf{N} = \left( \begin{array}{ccc}
1 & a        &b\\
c & n     &m\\
p & q       &r
\end{array} \right)
\end{displaymath}\\
D'apr\'es le th\'eor\`eme pr\'ec\`edant $Cat(N)\neq{\Oe}$ alors il
existe une cat\'egorie $A$ d\'efinie comme pr\'ec\`edemment, soit $A'$
une cat\'egorie dont les objets $Ob(A')=Ob(A)$ et les morphismes
$Mor(A')=Mor(A)\cup\{\lambda^0\lambda^0n^1,\lambda^0\lambda^0n^2,...,\lambda^0\lambda^0n^{z-1}\}$
alors
$A'(x_1,x_1)=\{id_{\lambda^0},\lambda^0\lambda^0n^1,\lambda^0\lambda^0n^2,...,\lambda^0\lambda^0n^{z-1}\}$
les \'equations de la loi
de composition associ\'ee \`a $n_i$ d\'efinies par:\\
$\lambda^0\lambda^0n^i(...)=(...)$\\
$(...)\lambda^0\lambda^0n^i=(...)$\\
$(\lambda^0\lambda^0n^i)(\lambda^0\lambda^0n^j)=(\lambda^0\lambda^0n^1)$avec $i\neq j$ \\
$(\lambda^0\lambda^0n^i)^2=(\lambda^0\lambda^0n^i)$ avec $i\neq j$ \\
donc $A'$ est une cat\'egorie associ\'e \`a M donc
$Cat(M)\neq{\Oe}$.\\
\\
\section {Matrices G\'en\'erales strictementes positives}
\begin{theoreme}
\label{gen1pos}
Soit $M$ une matrice de taille $n$ telle que $M$ d\'efinie par:\\
\begin{displaymath}
\mathbf{M} = \left( \begin{array}{cccc}
1 & M_{12} & \ldots &M_{1n}\\
M_{21} & M_{22} & \ldots & M_{2n}\\
\vdots & \vdots & \ddots&\vdots\\
M_{n1} & M_{n2} &  \ldots   &M_{nn}
\end{array} \right)
\end{displaymath}
avec $M_{ij}> 0$ $ \forall i,j \in\{1,...,n\}$ et $M_{ii}>1 $ pour $\in\{2,...,n\}$\\
alors $Cat(M)\neq {\Oe}$ si et seulement si $M_{ii}> M_{1i}M_{i1}
\forall i \in\{2,...,n\}$ et $M_{ij}\geq M_{i1}M_{1j}$ avec $i,j
\in\{2,...,n\}$\\
\end{theoreme}
$\Rightarrow)$ On suppose que $Cat(M)\neq {\Oe}$ alors il existe
$A$ une cat\'egorie associ\'ee \`a $M$ dont les objets sont
$\{\lambda^1,...,\lambda^n\}$.\\
On va d\'emontrer que $M_{ii}>M_{1i}M_{i1}$ on suppose que
$M_{ii}\leq M_{1i}M_{i1}$,soient\\
$A(\lambda^0,\lambda^i)=\{(\lambda^0\lambda^iG^{1,1}),...,(\lambda^0\lambda^iG^{1,a})\}$\\
$A(\lambda^i,\lambda^0)=\{(\lambda^i\lambda^0G^{1,1}),...,(\lambda^i\lambda^0G^{b,1})\}$\\
$A(\lambda^i,\lambda^i)=\{id_{\lambda^i},(\lambda^i\lambda^iG^1),...,(\lambda^i\lambda^iG^c)\}$\\
avec  $a=M_{1i}$ , $b=M_{i1}$ et $c=M_{ii}$ .\\
Ona
$(\lambda^i\lambda^0G^{i,1})(\lambda^0\lambda^iG^{1,j})$=$id_{\lambda^0}$
et
$(\lambda^0\lambda^iG^{1,j})(\lambda^i\lambda^0G^{i,1})=(\lambda^i\lambda^iG^p)$
Ona on a $M_{ii}\leq M_{1i}M_{i1}$  alors soit il existe $1\leq
i\neq j\leq b$,$1\leq k\leq a$ et $1\leq t\leq c$ tel que
$(\lambda^0\lambda^iG^{1,k})(\lambda^i\lambda^0G^{i,1})=
(\lambda^0\lambda^iG^{1,k})(\lambda^i\lambda^0G^{j,1})=
(\lambda^i\lambda^iG^t)$ alors\\
$(\lambda^i\lambda^0G^{i,1})\Big[(\lambda^0\lambda^iG^{1,k})(\lambda^i\lambda^0G^{i,1})\Big]$=
$(\lambda^i\lambda^0G^{i,1})\Big[(\lambda^0\lambda^iG^{1,k})(\lambda^i\lambda^0G^{j,1})\Big]$
donc\\
$\Big[(\lambda^i\lambda^0G^{i,1})(\lambda^0\lambda^iG^{1,k})\Big](\lambda^i\lambda^0G^{i,1})$=
$\Big[(\lambda^i\lambda^0G^{i,1})(\lambda^0\lambda^iG^{1,k})\Big](\lambda^i\lambda^0G^{j,1})$
alors\\ $(\lambda^i\lambda^0G^{i,1})=(\lambda^i\lambda^0G^{j,1})$,
donc  i=j introduction alors $M_{ii}>M_{1i}M_{i1}$.\\
Pour $M_{ij}\geq M_{i1}M_{1j}$\\
soient $A(\lambda^i,\lambda^0)=\{(\lambda^i\lambda^0G^{1,1}),...,(\lambda^i\lambda^0G^{b,1})\}$\\
$A(\lambda^0,\lambda^j)=\{(\lambda^1,\lambda^jG^{1,1}),...,(\lambda^1,\lambda^jG^{1,m})\}$\\
$A(\lambda^i,\lambda^j)=\{(\lambda^i\lambda^jG^1),...,(\lambda^i\lambda^jG^v)\}$\\
 avec $b=M_{i1},m=M_{1j}$ et
$v=M_{ij}$,on suppose que $M_{ij}< M_{i1}M_{1j}$  alors il existe
$1\leq i\neq j\leq b$,$1\leq k\neq c\leq m$ et $1\leq t\leq v$ tel
que:\\
$(\lambda^i\lambda^jG^t)\}=(\lambda^0,\lambda^jG^{1,k})(\lambda^i\lambda^0G^{i,1})=
(\lambda^0,\lambda^jG^{1,c})(\lambda^i\lambda^0G^{j,1})$\\
alors\\
$(\lambda^i\lambda^jG^t)(\lambda^0\lambda^iG^{1,a})=
\Big[(\lambda^0,\lambda^jG^{1,k})(\lambda^i\lambda^0G^{i,1})\Big](\lambda^0\lambda^iG^{1,a})=
\Big[(\lambda^0,\lambda^jG^{1,c})(\lambda^i\lambda^0G^{j,1})\Big](\lambda^0\lambda^iG^{1,a})$\\
donc\\
$(\lambda^i\lambda^jG^t)(\lambda^0\lambda^iG^{1,a})=
(\lambda^0,\lambda^jG^{1,k})\Big[(\lambda^i\lambda^0G^{i,1})(\lambda^0\lambda^iG^{1,a})\Big]=
(\lambda^0,\lambda^jG^{1,c})\Big[(\lambda^i\lambda^0G^{j,1})(\lambda^0\lambda^iG^{1,a})\Big]$\\
alors $(\lambda^0,\lambda^jG^{1,k})=(\lambda^0,\lambda^jG^{1,c})$
impossible car $k\neq c$ donc $M_{ij}\geq M_{i1}M_{1j}$ .\\
$\Leftarrow)$ soit M$'$ une matrice d'order n d\'efinie par:
\begin{displaymath}
\mathbf{M'} = \left( \begin{array}{cccc}
1 & M_{12} & \ldots &M_{1n}\\
M_{21} & (M_{21} M_{12}) & \ldots & (M_{21} M_{1n})\\
\vdots & \vdots & \ddots&\vdots\\
M_{n1} & (M_{n1}M_{12}) &  \ldots   &(M_{n1}M_{1n})
\end{array} \right)
\end{displaymath}
soit A$'$ une semi-cat\'egorie d\'efinie par:\\
$A'(\lambda^0,\lambda^0)$=$\{id_{\lambda^0}\}$\\
$A'(\lambda^i,\lambda^0)=\{(\lambda^i\lambda^0G^{1,1}),...,(\lambda^i\lambda^0G^{M_{i1},1})\}$\\
$A'(\lambda^0,\lambda^j)=\{(\lambda^1\lambda^jG^{1,1}),...,(\lambda^1,\lambda^jG^{1,M_{1i}})\}$\\
$A'(\lambda^i,\lambda^i)=\{(\lambda^i\lambda^iG^{1,1}),...,(\lambda^i\lambda^iG^{M_{1i},M_{i1}})\}$\\
$A'(\lambda^i,\lambda^j)=\{(\lambda^i\lambda^jG^{1,1}),...,(\lambda^i\lambda^jG^{M_{1j},M_{i1}})\}$ avec $i \neq j$\\
on d\'efinir la loi de composition par:\\
$(\lambda^j\lambda^kG^{u',v'})(\lambda^i\lambda^jG^{u,v})=(\lambda^i\lambda^jG^{u,v'})$\\
pour l'assiciativit\'e on a:\\
$\Big[(\lambda^j\lambda^kG^{u',v'})(\lambda^i\lambda^jG^{u'',v''})\Big](\lambda^k\lambda^pG^{u,v})=$\\
$(\lambda^i\lambda^kG^{u'',v'})(\lambda^k\lambda^pG^{u,v})=$\\
$(\lambda^k\lambda^kG^{u,v'})$\\
$(\lambda^j\lambda^kG^{u',v'})\Big[(\lambda^i\lambda^jG^{u'',v''})(\lambda^k\lambda^pG^{u,v})\Big]=$\\
$(\lambda^j\lambda^kG^{u',v'})(\lambda^k\lambda^jG^{u,v''})=$\\
$(\lambda^k\lambda^kG^{u,v'})$\\
Donc $\Big[(\lambda^j\lambda^kG^{u',v'})(\lambda^i\lambda^jG^{u'',v''})\Big](\lambda^k\lambda^pG^{u,v})=(\lambda^j\lambda^kG^{u',v'})\Big[(\lambda^i\lambda^jG^{u'',v''})(\lambda^k\lambda^pG^{u,v})\Big]$\\
Alors $A'$ semi-cat\'egorie associ\'ee \`a $M'$.\\
on ajoute des morphismes pour g\'en\`eraliser le th\'eor\`eme sur
les matrices de taille 3. On arrive surtout sur les  ensembles des morphismes suiventes:\\
$A''(\lambda^i,\lambda^i)=\{(\lambda^i\lambda^iG^{1,1}),...,(\lambda^i\lambda^iG^{M_{1i},M_{i1}}),(\lambda^i\lambda^iE^1),...,(\lambda^i\lambda^iE^{s_i})\}$  $ \forall i\in\{2,...,n\}$\\
$A''(\lambda^i,\lambda^j)=\{(\lambda^i\lambda^jG^{1,1}),...,(\lambda^i\lambda^jG^{M_{1j},M_{i1}}),(\lambda^i\lambda^jH^1),...,(\lambda^i\lambda^jH^{t^i_j})\}$$ \forall i \neq j$\\
\`a condition que tous les ajout\'es sont distincts, et en plus
la loi de composition est d\'efinie par :\\
$(\lambda^i\lambda^iE^k)\circ(...)=(\lambda^i\lambda^iG^{1,1})\circ (...)$ $\forall k\in\{1,...,s_i\}$ et $i\in\{2,...,n\}$\\
$(...)\circ(\lambda^i\lambda^iE^k)=(...) \circ (\lambda^i\lambda^iG^{M_{1i},M_{i1}}) $$\forall k\in\{1,...,s_i\}$ et $i\in\{2,...,n\}$\\
$(\lambda^i\lambda^jH^p)\circ(...)=(\lambda^i\lambda^jG^{1,1})\circ(...)$ $\forall p\in\{1,...,t^i_j\}$ et $i,j\in\{2,...,n\}$\\
$(...)\circ(\lambda^i\lambda^jH^p)=(...)\circ(\lambda^i\lambda^jG^{M_{1j},M_{i1}})$ $\forall p\in\{1,...,t^i_j\}$ et $i,j\in\{2,...,n\}$\\
$(\lambda^i\lambda^jH^p)\circ(\lambda^i\lambda^iE^k)=(\lambda^i\lambda^jG^{1,M_{i1}})$$\forall (p,k)\in\{1,...,t^i_j\}\times \{1,...,s_i\}$ et $i,j\in\{2,...,n\}$\\
$(\lambda^i\lambda^iE^k)\circ(\lambda^j\lambda^iH^p)=(\lambda^j\lambda^iG^{1,M_{i1}})$$\forall (p,k)\in\{1,...,t^j_i\}\times \{1,...,s_i\}$ et $i,j\in\{2,...,n\}$\\
Donc $A''$ une semi-cat\'egorie associ\'ee \`a $M''$ telle que cette matrice d\'efinie par:\\
\begin{displaymath}
\mathbf{M''} = \left( \begin{array}{cccc}
1 & M_{12} & \ldots &M_{1n}\\
M_{21} & (M_{21} M_{12})+s_i & \ldots & (M_{21} M_{1n})+t_n^2\\
\vdots & \vdots & \ddots&\vdots\\
M_{n1} & (M_{n1}M_{12}+)+t_2^n &  \ldots   &(M_{n1}M_{1n})+s_n
\end{array} \right)
\end{displaymath}
On peut ensuite ajouter les identit\'es sur $x_2,\ldots , x_n$. \\
\textbf{Finalement} si $M=(M_{ij})_n$ une matrice strictement
positive d'order $n$ telle que $M_{11}=1$et $m_ii>1$ alors
$Cat(M)\neq{\Oe}$ si et seulement si $M_{ii}> M_{1i}M_{i1} \forall
i \in\{1,...,n\}$ et $M_{ij}\geq
M_{i1}M_{1j} \forall i,j \in\{1,...,n\} i\neq j$.\\
\begin{corollary}: Pour toute matrice positive on peut \'etudier si
cette matrice marche ou non.
\end{corollary}
En effet: soit $M=(m_{ij})$ de taille n alors il y a deux cas :\\
a- $m_{ii}>1$ pour tout $ i \in\{1,2,...,n\}$\\
b- il existe au moins  une i' tel que $m_{i'i'}=1$\\
Cas (a):\\
on a $Cat(M)\neq{\Oe}$ d'apr\'es le th\'eor\`eme
pr\'ec\'edant.\\
Cas(b):\\
1- s'il existe une seule i' tel que $m_{i'i'}=1$ ,on peut
\'etudier
cette matrice d'apres le th\'eor\`eme 1.4.\\
2- s'ils existent i,j,......l tel que
$m_{ii}=m_{jj}=........=m_{ll}=1$ alors il y a deux cas:\\
-si M une matrice r\'eduite alors M ne pas marche d'apr\'es (lemme 4.6) .\\
-si M une matrice non r\'eduite alors il existe une matrice N
r\'eduction de M facile \`a \'etudier.\\
\\
\begin{corollary}
Si $M$ est une matrice de taille $n\geq 3$ avec $m_{ij}\geq 1$,
alors $Cat(M)\neq \Oe$ si et seulement si, pour toute sous-matrice
$N\subset M$ de taille $3$ on a $Cat(N)\neq \Oe$.
\end{corollary}
En effet, dans l'\'etude pr\'ec\'edente, dans les cas (a) et (b1)
les conditions ne concernent que les triples d'indices $i,j,k$ et
donc ne concernent que les sous-matrices de taille $3$. Pour cas
(b2) si $m_{ii}=m_{jj}=........=m_{ll}=1$, la condition
n\'ecessaire et suffisante pour que $M$ marche est que pour tout
autre $k$ on a $m_{ik}= m_{jk} = ... = m_{lk}$ et $m_{ki}= m_{kj}
= ... = m_{kl}$, et que la sous-matrice d\'efinie en enlevant
$j,...,l$ marche d'apr\`es le cas (b1).
\\

\section{Cas g\'en\'erale}
Soit $M$ une matrice d'ordre $n$ d\'efinie par:
\begin{displaymath}
\mathbf{M} = \left( \begin{array}{cccc}
m_{11} & m_{12} & \ldots &m_{1n}\\
m_{21} & m_{22} & \ldots & m_{2n} \\
\vdots & \vdots & \ddots&\vdots\\
m_{n1} & m_{n2} &  \ldots   &m_{nn}
\end{array} \right)
\end{displaymath}
avec $m_{ij}\geq 0 \forall (i,j)\in n\times n$ \\
On va chercher les classes et apr\`es d\'eterminer les
sous-matrices qui  d\'ecomposent la matrice pour isoler les
coificients nulles. Le r\'esultat de \cite{Allouch} s'applique \`a
chaque sous-matrice diagonale \`a coefficients strictement
positifs $M^{\lambda}$. Pour les sous-matrices $M^{\lambda , \mu}$
quand $\lambda > \mu$ la condition est tir\'ee du th\'eor\`eme
\ref{ex1} ci-dessous. On suppose $M$ une matrice qui est reduite au
sens de \cite{Allouch} telle que la relation d'ordre soit
transitive.\\\\
\textbf{Pour n=4}\\
Soit M une matrice d\'efinie par\\
\begin{displaymath} \mathbf{M} = \left(
\begin{array}{cc|cc}
1& b&c&d  \\
e&f&k&l\\
\hline
0&0&1&x\\
0&0&q&m
\end{array} \right) .
\end{displaymath}
avec $b,c,d,e,f,k,l,x,q$,et $m$  strictement positives.\\
soit $A=\{x_1,x_2,x_3,x_4\}$ associ\'e \`a  $M$, alors $Ob(A)/_\sim$=$\{1,2\}$ avec les classes $1=\{x_1,x_2\}$ et $2=\{x_3,x_4\}$.\\
Il y a 3 sous-matrices $M^{1}$,$M^{2}$ et $M^{1,2}$ qui d\'ecomposent la matrice $M$.\\
Dans les notations g\'en\'erales la matrice $M$ serait :
\begin{displaymath} \mathbf{M} = \left(
\begin{array}{cc|cc}
M(1^0, 1^0) & M(1^0, 1^1) & M(1^0, 2^0) & M(1^0, 2^1) \\
M(1^1, 1^0) & M(1^1, 1^1) & M(1^1, 2^0) & M(1^1, 2^1) \\
\hline
0 & 0 & M(2^0, 2^0) & M(2^0, 2^1) \\
0 & 0 & M(2^1, 2^0) & M(2^1, 2^1)
\end{array} \right) = \left(
\begin{array}{c|c}
\triangle_0 & \triangle^1_0 \\
\hline
 0&\triangle_1
\end{array} \right) .
\end{displaymath}
\begin{theoreme}
\label{ex1}
\begin{displaymath}  Cat(M)\neq{\emptyset}\Leftrightarrow\left\{
\begin{array}{llll}
f\geq be+1\\
m\geq xq+1\\
k,d,l\geq c \\
l\geq k \\
l\geq d \\
l \geq k+d-c
\end{array} \right.
\end{displaymath}
\end{theoreme}
Preuve:\\
 $\Rightarrow)$ Les deux sous-matrices positives
\begin{displaymath} \mathbf{M} = \left(
\begin{array}{cc}
1& b \\
e&f
\end{array} \right).
\end{displaymath}
et \begin{displaymath} \mathbf{M} = \left(
\begin{array}{cc}
1& x \\
q&m
\end{array} \right) .
\end{displaymath}
donnent les deux \'equations suivantes:\\
$f\geq be+1, m\geq xq+1$.\\
$Cat(M)\neq{\emptyset}$ alors d$'$apr\`es la m\'ethode
pr\'ec\`edente  il existe  $A$ une cate\'gorie donts les objets
sont $\{1^0,1^1,2^0,2^1\}$ avec $A/_\sim=\{1,2\}$ et les
morphismes
definis par:\\
Les morphismes associe\'es \`a $\triangle_0$:\\
$A(1^0,1^0)=\{id_{1^0}=1\}$\\
$A(1^0,1^1)=\{1^01^1G^{1,1},...,1^01^1G^{1,b}\}$\\
$A(1^1,1^0)=\{1^11^0G^{1,1},...,1^11^0G^{e,1}\}$\\
$A(1^1,1^1)=\{ id _{1^1}\}\cup
\{1^11^1G^{1,1},...,1^11^1G^{b,e}\}$
\\
$\cup \{1^11^1X^{1},...,1^11^1X^{f-be-1}\}$\\
\\
Les morphismes associe\'es \`a $\triangle_1$:\\
$A(2^0,2^0)=\{id_{2^0}=1=2^02^0G^{1,1}\}$\\
$A(2^0,2^1)=\{2^02^1G^{1,1},...,2^,2^1G^{1,x}\}$\\
$A(2^1,2^0)=\{2^12^0G^{1,1},...,2^12^0G^{q,1}\}$\\
$A(2^1,2^1)=\{ id _{2^1}\}\cup\{2^12^1G^{1,1},...,2^12^1G^{x,q}\}$
\\ $\cup \{2^12^1X^{1},...,2^12^1X^{m-xz-1}\}$.

Pour le moment, notons ainsi
les morphismes associe\'es \`a $\triangle_1^0$:\\
$A(1^0,2^0)=\{1^02^0A^1,...,1^02^0A^c\}$ notera par A\\
$A(1^1,2^0)=\{1^12^0B^1,...,1^12^0B^k\}$ notera par B\\
$A(1^0,2^1)=\{1^02^1C^1,...,1^02^1C^{d}\}$ notera par C\\
$A(1^1,2^1)=\{1^12^1D^1,...,1^12^1D^{l}\}$ notera par D
\\
On va demontrer que $k\geq c$.\\
on suppose que $k <c$ alors
$A(1^0,2^0)=\{1^02^0A^1,...,1^02^0A^k,...,1^02^0A^c\}$ et
$A(1^1,2^0)=\{1^12^0B^1,...,1^12^0B^k\}$ comme $k <c$ alors il
existe au moins $p\in \{1,...,k\}$, $u'\neq u'' \in \{1,...,c\}$ tel que:\\
$(1^01^1G^{1,1})(2^01^0A^{u'}) =(1^01^1G^{1,1})(2^01^0A^{u''})=(2^01^1B^p)\Rightarrow$\\
$(1^11^0G^{1,1})\Big[(1^01^1G^{1,1})(2^01^0A^{u'})\Big]=(1^11^0G^{1,1})\Big[(1^01^1G^{1,1})(2^01^0A^{u''})\Big]=(1^01^1G^{1,1})(2^11^0C^p)$
alors
$\Big[(1^11^0G^{1,1})(1^01^1G^{1,1})\Big](2^01^0A^{u'})=\Big[(1^11^0G^{1,1})(1^01^1G^{1,1})\Big](2^01^0A^{u''})=(1^01^1G^{1,1})(2^11^0C^p)$
donc $(2^01^0A^{u'})=(2^01^0A^{u''})=(1^01^1G^{1,1})(2^11^0C^p)$
alors $u'=u''$ contradiction donc  $k\geq c$. \\
la  meme pour les autres \'equations alors $d,l\geq c$.\\
Pour les autres in\'equations $l\geq k,d$\\
on va d\'emontrer que $l\geq k$.\\
On suppose que $l<k$ alors il existe au moins $p\in \{1,...,l\}$, $u'\neq u'' \in \{1,...,k\}$ tel que:\\
$(2^12^0G^{1,1})(2^01^1B^{u'}) =(2^12^0G^{1,1})(2^01^1B^{u''})=(2^11^1D^p)\Rightarrow$\\
$(2^02^1G^{1,1})\Big[(2^12^0G^{1,1})(2^01^1B^{u'})\Big]=(2^02^1G^{1,1})\Big[(2^12^0G^{1,1})(2^01^1B^{u''})\Big]=
(2^02^1G^{1,1})(2^11^1D^p)$ alors
\Big[$(2^02^1G^{1,1})(2^12^0G^{1,1})\Big](2^01^1B^{u'})=\Big[(2^12^0G^{1,1})(2^02^1G^{1,1})\Big](1^12^0B^{u''})$
donc $(1^12^0B^{u'})=(1^12^0B^{u''})$ alors $u'=u''$ contradiction
donc $l\geq k$ et aussi $l\geq d$.
Il reste l'\'equation $l\geq k+d-c$.\\
D'abord on va d\'emontrer que $(C-A)\subset D$ \\
et $(B-A)\subset D$ disjoints,\\
ce qui conduit a l'inegalite.\\
Pour voir qu'ils sont disjoints, fixons une meilleure notation:\\
on a:\\
$1^01^1G^{1,1}:1^0\rightarrow 1^1$, on note $1^01^1G^{1,1}$ par $G_1$ \\
$1^11^0G^{1,1}:1^1\rightarrow 1^0$, on note $1^11^0G^{1,1}$ par $F_1$ \\
Donc $(1^11^0G^{1,1})(1^01^1G^{1,1})=F_1G_1= id _{1^0}$\\
$2^02^1G^{1,1}:2^0\rightarrow 2^1$, on note $2^02^1G^{1,1}$ par $N_1$ \\
$2^12^0G^{1,1}:2^1\rightarrow 2^0$, on note $2^12^0G^{1,1}$ par $M_1$ \\
alors $(2^02^1G^{1,1})(2^12^0G^{1,1})=M_1N_1= id _{2^0}$\\
On d\'efinie les applications suivantes:\\
\xymatrix{gN_1:B  \ar[r]& D} tel que $gN_1(1^12^0B^i)=N_1(1^12^0B^i)=(2^02^1G^{1,1})(1^12^0B^i).$\\
\xymatrix{gM_1:D  \ar[r]& B} tel que $gM_1(1^12^1D^i)=M_1(1^12^1D^i)=(2^12^0G^{1,1})(1^12^1D^i).$\\
\xymatrix{dF_1:C  \ar[r]& D} tel que $dF_1(1^02^1C^i)=(1^02^1C^i)F_1=(1^02^1C^i)(1^11^0G^{1,1}).$\\
\xymatrix{dG_1:D  \ar[r]& C} tel que $dG_1(1^12^1D^i)=(1^12^1D^i)G_1=(1^12^1D^i)(1^01^1G^{1,1}).$\\
On remarque que les composes:\\
 \xymatrix{B  \ar[r]^{gN_1}& D
\ar[r]^{gM_1}&B}    $gM_1\circ gN_1=id_B$.\\
et\\
\xymatrix{D  \ar[r]^{dG_1}& C
\ar[r]^{dF_1}&D} $dF_1\circ dG_1=id_D$.\\
\begin{lemma}
\label{ldkc}:\\
Les applications $gN_1$ et $dF_1$ sont injectives, en plus
$gN_1(1^12^0B^i)=(2^02^1G^{1,1})(1^12^0B^i)=(1^12^1D^i)$ et
$dF_1(1^02^1C^i)=(1^02^1C^i)(1^11^0G^{1,1})=(1^12^1D^i)$
\end{lemma}
Preuve:\\
Soient $i\neq i'$ tel que $gN_1(1^12^0B^i)=gN_1(1^12^0B^{i'})$\\
$(2^02^1G^{1,1})(1^12^0B^i)=(2^02^1G^{1,1})(1^12^0B^{i'})$ donc \\
$\big[(2^12^0G^{1,1})(2^02^1G^{1,1})\big](1^12^0B^i)=\big[(2^02^1G^{1,1})(2^02^1G^{1,1})\big](1^12^0B^{i'})$
alors \\
$(1^12^0B^i)=(1^12^0B^{i'})$ donc $gN_1$ est injectife ce qui donne $gN_1(1^12^0B^i)=(2^02^1G^{1,1})(1^12^0B^i)=(1^12^1D^i)$ .\\
D'autre part:\\
Soient $i\neq i'$ tel que $dF_1(1^02^1C^i)=dF_1(1^02^1C^{i'})$\\
$(1^02^1C^i)(1^11^0G^{1,1})=(1^02^1C^{i'})(1^11^0G^{1,1})$ donc \\
$(1^02^1C^i)\big[(1^11^0G^{1,1})(1^01^1G^{1,1})\big]=(1^02^1C^{i'})\big[(1^11^0G^{1,1})(1^01^1G^{1,1})\big]$ alors \\
$(1^02^1C^i)=(1^02^1C^{i'})$ donc $dF_1$ est injectife ce qui donne $dF_1(1^02^1C^i)=(1^02^1C^i)(1^11^0G^{1,1})=(1^12^1D^i)$.\\
D'autre part on a le carre suivant:\\
$$ \xymatrix{
    A \ar[r]^{gN_{1}} \ar[d]_{dF_{1}} & C \ar[d]^{dF_{1}} \\
    B \ar[r]_{gN_{1}} & D
  }$$
o\`u les fleches verticales sont $dF_1$ et les fleches horizontales sont $gN_1$.\\
Le diagramme commute car  \\
$dF_1gN_1(u)=dF_1(N_1.u)= N_1uF_1= gN_1dF_1(u).$\\
Aussi les fleches sont tous injectives (comme ci-dessus).\\
On peut donc considerer $dF_1(gN_1(A))\subset D$. On notera cela par $N_1.A.F_1$.\\
On a egalement $dF_1(C)\subset D$, note par $C.F_1$,\\
similairement $gN_1(B)=N_1.B \subset D$.\\
On peut maintenant dire plus precisement ce qu'on veut demontrer:\\
que \\
($C.F_1 - N_1.A.F_1$) $\subset
D$ \\
et\\
($N_1.B - N_1.A.F_1$)$\subset D$ \\sont disjoints.\\
Par absurd, supposons qu'il existe $y\in(C.F_1 - N_1.A.F_1)\cap(N_1.B - N_1.A.F_1)$\\
$y\in (C.F_1 - N_1.A.F_1)$ ce qui donne  il existe $u\in (C - N_1.A)$ tel que $y=dF_1(u)=uF_1$ car $dF_1$ injectif\\
$y\in (N_1.B - N_1.A.F_1)$ ce qui donne  il existe $v\in (B -A.F_1)$ tel que $y=gN_1(v)=N_1v$ car $gN_1$ injectif\\
$y=uF_1 =N_1v$ \\
On considere alors \\
$gN_1.gM_1(u)=N_1M_1u =N_1M_1uF_1G_1=N_1M_1N_1vG_1=N_1vG_1=uF_1G_1= u$\\
On a $u\in (C - N_1.A)$ c.d.\`a $u\in C$ alors il existe a tel que
$u=(1^02^1C^a)$.\\
D'autre part $gM_1(u)=(2^12^0G^{1,1})(1^02^1C^a)=(1^02^0A^a)$ ce
qui donne $gM_1(u)\in A$ alors   $gN_1.gM_1(u)=u\in N_1.A$ ; une contradiction.\\
alternativement, de facon similaire \\
$dF_1.dG_1(v)=vG_1F_1 = M_1N_1vG_1F_1= M_1uF_1G_1F_1=M_1uF_1=M_1N_1v=v $\\
$v\in B$ alors il existe b tel que $v=(1^12^0B^b)$.\\
D'autre part $dG_1(v)=(1^12^0B^b)(1^01^1G^{1,1})=(1^02^0A^b)$ ce
qui donne $dG_1(v)\in A$
alors $dF_1.dG_1(v)=v\in A.F_1$, encore une contradiction.\\
Donc, ($C.F_1 - N_1.A.F_1$) $\subset
D$ \\
et\\
($N_1.B - N_1.A.F_1$) $\subset D$\\ sont disjoints.
\begin{lemma}
\label{ldkc2}:\\
Soit   \xymatrix{\relax \txt{A} \ar[r]^-{f} & B} une application
injectife tel que A,B deux ensembles finis alors:\\
Pour tout $X,Y\subset A $ alors $X\simeq f(X)$ ce qui donne
Card(X)= Card(f(X) et  $Card(X-Y)=Card(X)-Card(X \cap
Y)=Card(f(X))-Card(f(X \cap Y)).$
\end{lemma}
Preuve:\\
On a f une application injectife et $X\subset A$ alors la
\xymatrix{\relax \txt{X} \ar[r]^-{f/X} & f(X)} est bien bijectife
donc Card(X)=Card(f(X)) et $Card(X-Y)=Card(X)-Card(X \cap
Y)=Card(f(X))-Card(f(X \cap Y))$.\\
Comme $gN_1$ et $dF_1$ sont injectives alors d'apr\'es le lemme
pr\'ec\`edante alors:\\
$Card(N_1.A.F_1)=Card(gN_1(dF_1(A)))=Card(dF_1(A))=Card(A)=c.$\\
$Card((C.F_1)\cap(N_1.A.F_1))=Card(N_1.A.F_1)=c$ car
$(N_1.A.F_1)\subset(C.F_1)$.\\
$Card(C.F_1)=Card(dF_1(C))=Card(C)=d$ et $Card(N_1.B)=Card(gN_1(B))=Card(B)=k.$\\
D'autre part:\\
$Card(D - N_1.A.F_1)=Card(D)-Card((C.F_1)\cap
(N_1.A.F_1))=l-c$\\
$Card(C.F_1 - N_1.A.F_1)=Card(C.F_1)-Card((C.F_1)\cap
(N_1.A.F_1))=d-c$\\
$Card(N_1.B - N_1.A.F_1)=Card(N_1.B)-Card((N_1.B)\cap
(N_1.A.F_1))=k-c$\\
$C.F_1\subset D$ alors ($C.F_1 - N_1.A.F_1$) $\subset (D-N_1.A.F_1)$ \qquad (1)\\
et \\
$N_1.B\subset D$ alors ($N_1.B - N_1.A.F_1$) $\subset (D-N_1.A.F_1)$  \qquad (2) \\
Finalement (1)+(2) donne que Card($C.F_1 - N_1.A.F_1$)+
Card($N_1.B - N_1.A.F_1$)$\leq$ Card$(D-N_1.A.F_1)$ alors
$(k-c)+(d-c)\leq (l-c)$ donc d+k-c$\leq l$.\\
$\Leftarrow)$ On va voir la d\'emonstration en bas dans le cas
g\'en\'erale.\\
\begin{definition}:Soit
\label{defacceptable} M une matrice d$'$ordre n, on dit que M est
{\em acceptable} si la relation $\geq$ est reflexive et
transitive, o\`u
$\lambda^j\geq \mu^j$ si M$(\lambda^i,\mu^j)\geq 1$.\\
\end{definition}

Pour $\lambda \in U$ posons $a(\lambda ^i):= M(\lambda ^i,\lambda ^0)$
et $b(\lambda ^j):= M(\lambda ^0,\lambda ^j)$. 

\begin{theoreme}
\label{gen1}Soit M une matrice reduite d'order n alors
\begin{displaymath}
Cat(M)\neq{\emptyset}\Leftrightarrow \left\{
\begin{array}{llllll}
\qquad \qquad M  & \mbox{acceptable} &  \\
M(\lambda^{i},\lambda^{i})\geq & a(\lambda^{i}) b(\lambda^{i}) +1& \forall \lambda \in U , i\geq 1\\
M(\lambda^{i},\lambda^{j})\geq & a(\lambda^{i}) b(\lambda^{j}) & \forall \lambda \in U \\
M(\lambda^{i},\mu^{j}) \geq & M(\lambda^{i},\mu^{0}) & \forall \lambda > \mu , \mu \in U  \\
M(\lambda^{i},\mu^{j}) \geq & M(\lambda^{0},\mu^{j}) & \forall \lambda > \mu , \lambda \in U  \\
M(\lambda^{i},\mu^{j}) \geq  & M(\lambda^{0},\mu^{j})+
M(\lambda^{i},\mu^{0})
& \\
& \qquad -M(\lambda^{0},\mu^{0}) &\forall \lambda \geq \mu \in U
\end{array} \right.
\end{displaymath}
\end{theoreme}
Preuve:\\
$\Rightarrow)$:$Cat(M)\neq\emptyset$ on a d'apr\`es la partition
de $M$ alors il existe A une cat\'egorie associ\'e \`a M dont les
objets
sont $\{\lambda^i,.....,\mu^j,.....,........etc \}$ \\
Donc M sera:
\begin{displaymath}
M=\bigcup_{\substack{\lambda\in U }}M^{\lambda}\bigoplus
\bigcup_{\substack{\mu\in V}}M^{\mu}\bigoplus
\bigcup_{\begin{subarray}{l}\lambda\in U\\
\mu\in V
\end{subarray}}M^{\lambda,\mu}
\end{displaymath}
On va d\'emontrer que $M$ est acceptable.\\
En effet soit $\lambda^i\in A$ comme $M^{\lambda}>0$ alors $M(\lambda^i,\lambda^i)\geq 1$ donc $\lambda^i\geq \lambda^i$ alors $\geq$ est reflexive.\\
Pour la transitivit\'e soient $\lambda^i$,$\mu^j$ et $\phi^k$ trois objets telque $\lambda^i \geq \mu^j$ et $\mu^j\geq \phi^k$ .\\
alors il existe deux morphismes F=$(\lambda^i,\mu^j,H^a)$ et
G=$(\mu^j,\phi^k,J^b)$. \\
On a
K=$(\mu^j,\phi^k,J^b)(\lambda^i,\mu^j,H^a)=(\lambda^i,\phi^k,M^c)$
donc $\lambda^i\geq \phi^k$ alors $\geq$ est transitive.\\
D'autre part le r\'esultat de \cite{Allouch2} s'applique \`a
chaque sous-matrice diagonale \`a coefficients strictement
positifs $M^{\lambda}$,$M^{\mu}$. Pour les sous-matrices
$M^{\lambda , \mu}$ quand $\lambda > \mu$ la condition est tir\'ee
du
th\'eor\`eme \ref{ex1} ci-dessus.\\
$ \Leftarrow)$:nous avons les conditions on va d\'emontrer
$Cat(M)\neq\emptyset$.\\
D'apres la m\'etode d'ajouter des elements voir \cite{Allouch}
nous pouvons travailler sur $\bigsqcup$ c.\`a.d M sera:
\begin{displaymath} M=\bigcup_{\substack{\lambda\in U \bigcup
V}}M^{\lambda}\bigoplus \bigcup_{\begin{subarray}{l}\lambda,\mu\in
U \bigcup V
\end{subarray}}M^{\lambda,\mu}
\end{displaymath}
Comme dans \cite{Allouch2}, il suffit de consid\'erer la
construction d'une {\em semi-cat\'egorie partiellement unitaire}
avec condition d'unit\'e seulement sur les $A(\lambda ^0, \lambda
^0)$ pour $\lambda \in U$, et en supposant que $M(\lambda
^i,\lambda ^j)=a(\lambda ^i)b(\lambda ^j)$ pour $\lambda \in U$.
Pour $\lambda \in V$ on peut prendre $M(\lambda ^i,\lambda ^j)=1$.

Soit A une (semi-)cat\'egorie donts les objets sont $\{\lambda^i$
avec $\lambda \in A/\sim\}$ et les morphismes sont
$\{(\lambda^i,\lambda^j,(u,v))\}\cup\{(\lambda^i,\mu^j,k)\}$.\\

Ici $1\leq u\leq a(\lambda ^i)$ et $1\leq v\leq b(\lambda ^j)$,
en mettant $a(\lambda ^i)=b(\lambda ^j)=1$ si $\lambda \in V$.\\

Pour $\lambda >\mu$, on notera
$$
|\mu^0,\lambda^0| := \{ 1,2,\ldots , M(\mu ^0,\lambda ^0)\},
$$
$$
|\mu^0,\lambda^i| := \{ 1,2,\ldots , M(\mu ^0,\lambda ^i)\},
$$
mais en revanche
$$
|\mu^n,\lambda^0| := \{ 1+M(\mu ^0,\lambda ^0)-M(\mu ^n,\lambda
^0), \ldots , 0,1,2,\ldots , M(\mu ^0,\lambda ^0) \},
$$
et enfin
$$
|\mu^n,\lambda^i| := \{ 1+M(\mu ^0,\lambda ^0)-M(\mu ^n,\lambda
^0), \ldots , 0,1,2,\ldots , h \},
$$
o\`u $h=M(\mu ^n, \lambda ^i)+ M(\mu ^0,\lambda ^0)-M(\mu
^n,\lambda ^0)$. De cette fa\c{c}on,
$$
\# |\mu^0,\lambda^0| = M(\mu^0,\lambda^0),
$$
$$
\# |\mu^0,\lambda^i| = M(\mu^0,\lambda^i),
$$
$$
\# |\mu^n,\lambda^0| = M(\mu^n,\lambda^0),
$$
et
$$
\# |\mu^n,\lambda^i| = M(\mu^n,\lambda^i).
$$
D'autre part
$$
|\mu^0,\lambda^i| \cap |\mu^n,\lambda^0| = |\mu^0,\lambda^0|.
$$
On notera souvent par $H'(\mu^n,\lambda^i)$ le compl\'ementaire
$$
H'(\mu^n,\lambda^i):= \{ M(\mu ^0,\lambda ^i)+1, \ldots , h\}
$$
ayant
$$
\# H'(\mu^n,\lambda^i) = M(\mu^n,\lambda^i) - M(\mu^0,\lambda^i)-
M(\mu^n,\lambda^0) + M(\mu^0,\lambda^0)
$$
\'el\'ements. L'hypoth\`ese dit que ce nombre d'\'el\'ements est
$\geq 0$.

Ainsi $|\mu^n,\lambda^i|$ est la r\'eunion disjointe des quatres
parties suivantes:
$$
|\mu^0,\lambda^0|
$$
$$
|\mu^0,\lambda^i| -  |\mu^0,\lambda^0|
$$
$$
|\mu^n,\lambda^0| -  |\mu^0,\lambda^0|
$$
$$
H'(\mu^n,\lambda^i).
$$

La loi de composition interne sera d\'efinie par:\\
a\big)($\lambda^{i}$,$\mu^{j}$,k)$\circ$($\varphi^{n}$,$\lambda^{i}$,$k'$)=($\varphi^{n}$,$\mu^{j}$,1)\\
(l'utilisation des symboles diff\'erents $\varphi , \lambda , \mu
$ ici veut dire qu'on a l'hypoth\`ese $\varphi \neq \lambda \neq
\mu$, une convention en rigeur partout),
\\
\\
b\big)($\lambda^{i}$,$\lambda^{j}$,(u,v))$\circ$($\lambda^{n}$,$\lambda^{i}$,($u',v'$))=($\lambda^{n}$,$\lambda^{j}$,$(u',v$))\\
\\
c\big)$(\lambda^i,\lambda^j,(u,v))\circ(\mu^n,\lambda^i,k)$=$(\mu^n,\lambda^j,k')$\\
Si $\lambda \in V$ alors $(\lambda^i,\lambda^j,(u,v))\circ(\mu^n,\lambda^i,k)$=$(\mu^n,\lambda^j,1)$\\
Si $\lambda \in U$ alors on a deux  cas :\\
Si i=0 alors $k'=k$ \qquad \\
Si i$\neq 0$ alors
\begin{displaymath}
k' = \left\{
\begin{array}{ll}
k & \textrm{ si $ k\in |\mu^0,\lambda^0|$}\\
1 & \textrm{ si $ k\in |\mu^0,\lambda^i|-|\mu^0,\lambda^0|$}\\
k & \textrm{ si $ k\in |\mu^n,\lambda^0|-|\mu^0,\lambda^0|$}\\
1 & \textrm{ si $ k\in H'(\mu^n,\lambda^i )$}
\end{array} \right.
\end{displaymath}
\\
d\big)$(\lambda^j,\mu^n,k)\circ(\lambda^i,\lambda^j,(u,v))$=$(\lambda^i,\mu^n,k')$\\
Si $\lambda \in V$ alors $(\lambda^j,\mu^n,k)\circ(\lambda^i,\lambda^j,(u,v))$=$(\lambda^i,\mu^n,1)$\\
Si $\lambda \in U$ on a deux  cas :\\
Si j=0 alors $k'=k$\\
Si $j\neq 0$ alors
\begin{displaymath}
k' = \left\{
\begin{array}{ll}
k & \textrm{ si $ k\in |\lambda^0,\mu^0|$}\\
1 & \textrm{ si $ k\in |\lambda^j,\mu^0|-|\lambda^0,\mu^0|$}\\
k & \textrm{ si $ k\in |\lambda^0,\mu^n|-|\lambda^0,\mu^0|$}\\
1 & \textrm{ si $ k\in H'(\mu^n,\lambda^i )$}
\end{array} \right.
\end{displaymath}
\\
\underline{Pour l'unit\'e partiel}:\\
Comme $k'=k$ pour $i=0$ (cas (c)) ou $j=0$ (cas (d)) quand
$\lambda \in U$, $(\lambda ^i \lambda ^i (1,1))$ agit comme
l'identit\'e.
\\
\underline{Pour l'associativit\'e}:\\
soient $\lambda$,$\mu$,$\varphi$,$\phi$ $\in A/\sim$ alors il y a
8 formes des \'equations associatives :\\
1)\quad$\lambda^i\mapsto\mu^j\mapsto\phi^t\mapsto\varphi^n$\\
2)\quad$\lambda^i\mapsto\mu^j\mapsto\mu^t\mapsto\phi^n$\\
3)\quad$\lambda^i\mapsto\mu^j\mapsto\phi^t\mapsto\phi^n$\\
4)\quad$\lambda^i\mapsto\lambda^j\mapsto\mu^t\mapsto\phi^n$\\
5)\quad$\lambda^i\mapsto\lambda^j\mapsto\mu^t\mapsto\mu^n$\\
6)\quad$\lambda^i\mapsto\lambda^j\mapsto\lambda^t\mapsto\mu^n$\\
7)\quad$\lambda^i\mapsto\mu^j\mapsto\mu^t\mapsto\mu^n$\\
8)\quad$\lambda^i\mapsto\lambda^j\mapsto\lambda^t\mapsto\lambda^n$\\
Rappelons que les symboles $\lambda , \mu , \phi , \varphi$
distincts repr\'esentent par convention des classes diff\'erents.
\\
On va v\'erifier les \'egalit\'es d'associativit\'e pour chacun de
ces cas.
\\
\underline{\'equation 1}\\
$(\phi^t,\varphi^n,a)(\mu^j,\phi^t,b)(\lambda^i,\mu^j,c)$\\
=$[(\phi^t,\varphi^n,a)(\mu^j,\phi^t,b)](\lambda^i,\mu^j,c)$\\
=$(\mu^j,\varphi^n,1)(\lambda^i,\mu^j,c)$\\
=$(\lambda^i,\varphi^n,1)$\\
=$(\phi^t,\varphi^n,a)[(\mu^j,\phi^t,b)(\lambda^i,\mu^j,c)]$\\
\underline{\'equation 2}\\
$\lambda^i\mapsto\mu^j\mapsto\mu^t\mapsto\phi^n$\\
Pour l$'$identit\'e:\\
$(\mu^0,\phi^n,a)(\mu^0,\mu^0,(1,1))(\lambda^i,\mu^0,b)$\\
=$[(\mu^0,\phi^n,a)(\mu^0,\mu^0,(1,1))](\lambda^i,\mu^0,b)$\\
=$(\mu^0,\phi^n,a)(\lambda^i,\mu^j,b)$\\
=$(\lambda^i,\phi^n,1)$\\
=$(\mu^0,\phi^n,a)[(\mu^0,\mu^0,(1,1))(\lambda^i,\mu^0,b)]$\\
=$(\mu^0,\phi^n,a)(\lambda^i,\mu^0,b)$\\
=$(\lambda^i,\phi^n,1)$\\
Si $\mu \in U$ alors:\\
$(\mu^t,\phi^n,a)(\mu^j,\mu^t,(u,v))(\lambda^i,\mu^j,b)$\\
=$[(\mu^t,\phi^n,a)(\mu^j,\mu^t,(u,v))](\lambda^i,\mu^j,b)$\\
=$(\mu^j,\phi^n,c)(\lambda^i,\mu^j,b)$\\
=$(\lambda^i,\phi^n,1)$\\
=$(\mu^t,\phi^n,a)[(\mu^j,\mu^t,(u,v))(\lambda^i,\mu^j,b)]$ facile \`a d\'emontrer\\
Si $\mu \in V$ alors:\\
$(\mu^t,\phi^n,a)(\mu^j,\mu^t,(u,v))(\lambda^i,\mu^j,b)$\\
=$[(\mu^t,\phi^n,a)(\mu^j,\mu^t,(u,v))](\lambda^i,\mu^j,b)$\\
=$(\mu^j,\phi^n,1)(\lambda^i,\mu^j,b)$\\
=$(\lambda^i,\phi^n,1)$\\
=$(\mu^t,\phi^n,a)[(\mu^j,\mu^t,(u,v))(\lambda^i,\mu^j,b)]$ \\
=$(\phi^t,\varphi^n,a)(\lambda^i,\mu^t,1)$\\
=$(\lambda^i,\phi^n,1)$\\
\underline{\'equation 3}\\
$(\phi^t,\phi^n,(u,v))(\mu^j,\phi^t,a)(\lambda^i,\mu^j,b)$\\
Pour l$'$identit\'e:\\
=[$(\phi^0,\phi^0,(1,1))(\mu^j,\phi^0,a)](\lambda^i,\mu^j,b)$\\
=$(\mu^j,\phi^0,c)(\lambda^i,\mu^j,b)$\\
=$(\lambda^i,\phi^0,1)$\\
=$(\phi^0,\phi^0,(1,1))[(\mu^j,\phi^0,a)(\lambda^i,\mu^j,b)]$\\
=$(\phi^0,\phi^0,(1,1))(\lambda^i,\phi^0,1)$\\
=$(\lambda^i,\phi^0,1)$\\
Si $\phi \in U$ alors: \\
=[$(\phi^t,\phi^n,(u,v))(\mu^j,\phi^t,a)](\lambda^i,\mu^j,b)$\\
=$(\mu^j,\phi^n,k)(\lambda^i,\mu^j,b)$\\
=$(\lambda^j,\phi^n,1)$\\
=$(\phi^t,\phi^n,(u,v))[(\mu^j,\phi^t,a)(\lambda^i,\mu^j,b)]$\\
=$(\phi^t,\phi^n,(u,v))(\lambda^i,\phi^t,1)$\\
=$(\lambda^i,\phi^n,1)$\\
Si $\phi \in V$ alors: \\
=[$(\phi^t,\phi^n,(u,v))(\mu^j,\phi^t,a)](\lambda^i,\mu^j,b)$\\
=$(\mu^j,\phi^n,1)(\lambda^i,\mu^j,b)$\\
=$(\lambda^j,\phi^n,1)$\\
=$(\phi^t,\phi^n,(u,v))[(\mu^j,\phi^t,a)(\lambda^i,\mu^j,b)]$\\
=$(\phi^t,\phi^n,(u,v))(\lambda^i,\phi^t,1)$\\
=$(\lambda^i,\phi^n,1)$\\
\underline{\'equation }4 comme l'\'equation 3\\
\underline{\'equation 5}\\
$(\mu^t,\mu^n,(a,b))(\lambda^j,\mu^t,k)(\lambda^i,\lambda^j,(c,d))$\\
Il y a 4 cas sur:\\
1)$\mu \in U$ et $\lambda \in V$\\
2)$\mu \in V$ et $\lambda \in U$\\
3)$\mu \in V$ et $\lambda \in V$\\
4)$\mu \in U$ et $\lambda \in U$\\
Cas 1) alors:\\
Pour l$'$identit\'e:\\
Q=$(\mu^0,\mu^0,(1,1))\big[(\lambda^j,\mu^0,k)(\lambda^i,\lambda^j,(c,d))\big]$\\
=$(\mu^0,\mu^0,(1,1))(\lambda^i,\mu^0,1)$\quad car $\lambda \in V$\\
=$(\lambda^i,,\mu^0,1)$\\
$Q'=\big[(\mu^0,\mu^0,(1,1))(\lambda^j,\mu^0,k)\big](\lambda^i,\lambda^j,(c,d))$\\
=$(\lambda^j,,\mu^0,k')(\lambda^i,\lambda^j,(c,d))$\\
=$(\lambda^i,,\mu^0,1)$\quad car $\lambda \in V$\\
Donc $Q=Q'$\\
$\big[(\mu^0,\mu^0,(1,1))(\lambda^0,\mu^0,k)\big](\lambda^0,\lambda^0,(1,1))=(\mu^0,\mu^0,(1,1))\big[(\lambda^0,\mu^0,k)(\lambda^0,\lambda^0,(1,1))\big]$
facile \`a d\'emontrer \\
Q=$(\mu^t,\mu^n,(a,b))\big[(\lambda^0,\mu^t,k)(\lambda^0,\lambda^0,(1,1))\big]$\\
=$(\mu^t,\mu^n,(a,b))(\lambda^i,\mu^t,1)$\quad car $\lambda \in V$\\
=$(\lambda^i,,\mu^n,1)$\\
$Q'=\big[(\mu^t,\mu^n,(a,b))(\lambda^j,\mu^t,k)\big](\lambda^i,\lambda^j,(c,d))$\\
=$(\lambda^j,,\mu^n,k')(\lambda^i,\lambda^j,(c,d))$\\
=$(\lambda^i,,\mu^n,1)$\quad car $\lambda \in V$\\
Donc $Q=Q'$\\
Q=$(\mu^t,\mu^n,(a,b))\big[(\lambda^j,\mu^t,k)(\lambda^i,\lambda^j,(c,d))\big]$\\
=$(\mu^t,\mu^n,(a,b))(\lambda^i,\mu^t,1)$\quad car $\lambda \in V$\\
=$(\lambda^i,,\mu^n,1)$\\
$Q'=\big[(\mu^t,\mu^n,(a,b))(\lambda^j,\mu^t,k)\big](\lambda^i,\lambda^j,(c,d))$\\
=$(\lambda^j,,\mu^n,k')(\lambda^i,\lambda^j,(c,d))$\\
=$(\lambda^i,,\mu^n,1)$\quad car $\lambda \in V$\\
Donc $Q=Q'$\\
Cas 2) comme 1)\\
Cas 3)\\
Pour l$'$identit\'e:
Q=$(\mu^0,\mu^0,(1,1))\big[(\lambda^j,\mu^0,k)(\lambda^i,\lambda^j,(c,d))\big]$\\
=$(\mu^0,\mu^0,(1,1))(\lambda^i,\mu^0,1)$\quad car $\lambda \in V$\\
=$(\lambda^i,\mu^0,1)$\\
$Q'=\big[(\mu^0,\mu^0,(1,1))(\lambda^j,\mu^0,k)\big](\lambda^i,\lambda^j,(c,d))$\\
=$(\lambda^j,,\mu^0,1)(\lambda^i,\lambda^j,(c,d))$ \quad car $\mu \in V$\\
=$(\lambda^i,,\mu^0,1)$\\
Alors $Q=Q'$\\
l$'$autre cas de l$'$identit\'e est le m\^eme.\\
Q=$(\mu^t,\mu^n,(a,b))\big[(\lambda^j,\mu^t,k)(\lambda^i,\lambda^j,(c,d))\big]$\\
=$(\mu^t,\mu^n,(a,b))(\lambda^i,\mu^t,1)$\quad car $\lambda \in V$\\
=$(\lambda^i,\mu^n,1)$\\
$Q'=\big[(\mu^t,\mu^n,(a,b))(\lambda^j,\mu^t,k)\big](\lambda^i,\lambda^j,(c,d))$\\
=$(\lambda^j,,\mu^n,1)(\lambda^i,\lambda^j,(c,d))$ \quad car $\mu \in V$\\
=$(\lambda^i,,\mu^n,1)$\\
Donc $Q=Q'$\\
Cas 4)\\
Q=$(\mu^t,\mu^n,(a,b))\big[(\lambda^0,\mu^t,k)(\lambda^0,\lambda^0,(1,1))\big]$\\
On a deux cas :\\
Si t=0\\
Q=$(\mu^0,\mu^n,(a,b))\big[(\lambda^0,\mu^0,k)(\lambda^0,\lambda^0,(1,1))\big]$\\
=$(\mu^0,\mu^n,(a,b))(\lambda^0,\mu^0,k)$\\
=$(\lambda^0,\mu^n,k)$\\
Q'=$\big[(\mu^0,\mu^n,(a,b))(\lambda^0,\mu^0,k)\big](\lambda^0,\lambda^0,(1,1))$\\
=$(\lambda^0,\mu^0,k)(\lambda^0,\lambda^0,(1,1))$\\
=$(\lambda^0,\mu^n,k)$\\
Alors $Q=Q'$\\
Si t$\neq 0$\\
Q=$(\mu^t,\mu^n,(a,b))\big[(\lambda^0,\mu^t,k)(\lambda^0,\lambda^0,(1,1))\big]$\\
=$(\mu^t,\mu^n,(a,b))(\lambda^0,\mu^t,k)$\\
=$(\lambda^0,\mu^n,k')$\\
\begin{displaymath}
avec \qquad k' = \left\{ \begin{array}{ll}
k & \textrm{ $|\lambda^0,\mu^0|$}\\
1 & \textrm{ $|\lambda^0,\mu^0|-|\lambda^0,\mu^0|$}\\
k & \textrm{ $|\lambda^0,\mu^t|-|\lambda^0,\mu^0|$}\\
1 & \textrm{ $H'(\lambda^0,\mu^t)$}
\end{array} \right.
\end{displaymath}
Donc k$'$=k sur  $|\lambda^0,\mu^t|$.\\
D$'$autre part:\\
Q$'=\big[(\mu^t,\mu^n,(a,b))(\lambda^0,\mu^t,k)\big](\lambda^0,\lambda^0,(1,1))$\\
=$(\lambda^0,\mu^n,k')(\lambda^0,\lambda^0,(1,1))$\\
=$(\lambda^0,\mu^n,k')$\\
\begin{displaymath}
avec \qquad k' = \left\{ \begin{array}{ll}
k & \textrm{ $|\lambda^0,\mu^0|$}\\
1 & \textrm{ $|\lambda^0,\mu^0|-|\lambda^0,\mu^0|$}\\
k & \textrm{ $|\lambda^0,\mu^t|-|\lambda^0,\mu^0|$}\\
1 & \textrm{ $H'(\lambda^0,\mu^t)$}
\end{array} \right.
\end{displaymath}
Alors $Q=Q'$\\
Pour
Q=$(\mu^0,\mu^0,(1,1))(\lambda^0,\mu^0,k)(\lambda^0,\lambda^0,(1,1))$
facile \`a d\'emontrer\\
Q=$(\mu^t,\mu^n,(a,b))\big[(\lambda^j,\mu^t,k)(\lambda^i,\lambda^j,(c,d))\big]$\\
et\\
$Q'=\big[(\mu^t,\mu^n,(a,b))(\lambda^j,\mu^t,k)\big](\lambda^i,\lambda^j,(c,d))$\\
On va v\'erifier $Q=Q'$ sur $|\lambda^j,\mu^t|$\\
Q=$(\mu^t,\mu^n,(a,b))\big[(\lambda^j,\mu^t,k)(\lambda^i,\lambda^j,(c,d))\big]$\\
=$(\mu^t,\mu^n,(a,b))(\lambda^i,\mu^t,k')$\\
=$(\lambda^i,\mu^n,k'')$\\
On a 4 possibilites :\\
1)j=0 et t=0\\
2)j=0 et t$\neq 0$\\
3)$j\neq 0$ et t=0\\
4)$j\neq 0$ et t$\neq 0$\\
$$j=0 , t=0$$
Q=$(\mu^0,\mu^n,(a,b))\big[(\lambda^0,\mu^0,k)(\lambda^i,\lambda^0,(c,d))\big]$\\
=$(\mu^0,\mu^n,(a,b))(\lambda^i,\mu^0,k)$ \quad voir  d) \\
=$(\lambda^i,\mu^n,k)$ \quad voir  c)\\
D'autre part\\
$Q'=\big[(\mu^0,\mu^n,(a,b))(\lambda^0,\mu^0,k)\big](\lambda^i,\lambda^0,(c,d))$\\
=$(\lambda^0,\mu^n,k)(\lambda^i,\lambda^0,(c,d))$ \quad voir  c)\\
=$(\lambda^i,\mu^n,k)$ \quad voir  d)\\
Donc $Q'=Q$ sur $|\lambda^0,\mu^0|$\\
$$j=0 , t\neq 0$$
Q=$(\mu^t,\mu^n,(a,b))\big[(\lambda^0,\mu^t,k)(\lambda^i,\lambda^0,(c,d))\big]$\\
=$(\mu^t,\mu^n,(a,b))(\lambda^i,\mu^t,k)$ \quad voir d)\\
=$(\lambda^i,\mu^n,k'')$ \quad voir c)\\
\begin{displaymath}
avec \qquad k'' = \left\{ \begin{array}{ll}
k & \textrm{ $|\lambda^0,\mu^0|$}\\
1 & \textrm{ $|\lambda^0,\mu^t|-|\lambda^0,\mu^0|$}\\
k & \textrm{ $|\lambda^i,\mu^0|-|\lambda^0,\mu^0|$}\\
1 & \textrm{ $H'(\lambda^i,\mu^t)$}
\end{array} \right.
\end{displaymath}
D'autre part,\\
$Q'=\big[(\mu^t,\mu^n,(a,b))(\lambda^0,\mu^t,k)\big](\lambda^i,\lambda^0,(c,d))$\\
=$(\lambda^0,\mu^n,k')(\lambda^i,\lambda^0,(c,d))$ voir c)\\
=$(\lambda^i,\mu^n,k')$\quad voir d)\\
\begin{displaymath}
avec \qquad k' = \left\{ \begin{array}{ll}
k & \textrm{ $|\lambda^0,\mu^0|$}\\
1 & \textrm{ $|\lambda^0,\mu^t|-|\lambda^0,\mu^0|$}\\
k & \textrm{ $|\lambda^0,\mu^0|-|\lambda^0,\mu^0|$}\\
1 & \textrm{ $H'(\lambda^0,\mu^t)$}
\end{array} \right. \quad = \quad\left\{ \begin{array}{ll}
k & \textrm{ $|\lambda^0,\mu^0|$}\\
1 & \textrm{ $|\lambda^0,\mu^t|-|\lambda^0,\mu^0|$}
\end{array} \right.
\end{displaymath}
Donc $Q=Q'$ sur $|\lambda^0,\mu^t|$\\
$$ j\neq 0, t=0$$
Q=$(\mu^0,\mu^n,(a,b))\big[(\lambda^j,\mu^0,k)(\lambda^i,\lambda^j,(c,d))\big]$\\
=$(\mu^0,\mu^n,(a,b))(\lambda^i,\mu^0,k')$ \quad voir d) \\
=$(\lambda^i,\mu^n,k')$ \quad voir c)\\
\begin{displaymath}
avec \qquad k' = \left\{ \begin{array}{ll}
k & \textrm{ $|\lambda^0,\mu^0|$}\\
1 & \textrm{ $|\lambda^j,\mu^0|-|\lambda^0,\mu^0|$}\\
k & \textrm{ $|\lambda^0,\mu^0|-|\lambda^0,\mu^0|$}\\
1 & \textrm{ $H'(\lambda^j,\mu^0)$}
\end{array} \right. \quad = \quad\left\{ \begin{array}{ll}
k & \textrm{ $|\lambda^0,\mu^0|$}\\
1 & \textrm{ $|\lambda^j,\mu^0|-|\lambda^0,\mu^0|$}
\end{array} \right.
\end{displaymath}
L'autre sens,\\
$Q'=\big[(\mu^0,\mu^n,(a,b))(\lambda^j,\mu^0,k)\big](\lambda^i,\lambda^j,(c,d))$\\
=$(\lambda^j,\mu^n,k)(\lambda^i,\lambda^j,(c,d))$ voir c)\\
=$(\lambda^i,\mu^n,k'')$ \quad voir d)\\
\begin{displaymath}
avec \qquad k'' = \left\{ \begin{array}{ll}
k & \textrm{ $|\lambda^0,\mu^0|$}\\
1 & \textrm{ $|\lambda^j,\mu^0|-|\lambda^0,\mu^0|$}\\
k & \textrm{ $|\lambda^0,\mu^n|-|\lambda^0,\mu^0|$}\\
1 & \textrm{ $H'(\lambda^j,\mu^n)$}
\end{array} \right.
\end{displaymath}
Alors $Q=Q'$ sur $|\lambda^j,\mu^0|$\\
$$ j\neq 0, t\neq 0$$
Q=$(\mu^t,\mu^n,(a,b))\big[(\lambda^j,\mu^t,k)(\lambda^i,\lambda^j,(c,d))\big]$\\
=$(\mu^t,\mu^n,(a,b))(\lambda^i,\mu^t,k')$ \quad voir d)\\
=$(\lambda^i,\mu^n,k'')$ \quad voir c)\\
\begin{displaymath}
avec \qquad k' = \left\{ \begin{array}{ll}
k & \textrm{ $|\lambda^0,\mu^0|$}\\
1 & \textrm{ $|\lambda^j,\mu^0|-|\lambda^0,\mu^0|$}\\
k & \textrm{ $|\lambda^0,\mu^t|-|\lambda^0,\mu^0|$}\\
1 & \textrm{ $H'(\lambda^j,\mu^t)$}
\end{array} \right.
\end{displaymath}
et
\begin{displaymath}
avec \qquad k'' = \left\{ \begin{array}{ll}
k' & \textrm{ $|\lambda^0,\mu^0|$}\\
1 & \textrm{ $|\lambda^0,\mu^t|-|\lambda^0,\mu^0|$}\\
k' & \textrm{ $|\lambda^i,\mu^0|-|\lambda^0,\mu^0|$}\\
1 & \textrm{ $H'(\lambda^i,\mu^t)$}
\end{array} \right. \quad = \quad\left\{ \begin{array}{ll}
k & \textrm{ $|\lambda^0,\mu^0|$}\\
1 & \textrm{ $|\lambda^0,\mu^t|-|\lambda^0,\mu^0|$}\\
1 & \textrm{ $|\lambda^j,\mu^0|-|\lambda^0,\mu^0|$}\\
1 & \textrm{ $H'(\lambda^i,\mu^t)$}
\end{array} \right.
\end{displaymath}
On va verifier que k$''$=1 sur
$|\lambda^j,\mu^0|-|\lambda^0,\mu^0|$.\\
logiquement $k''=k'$ ou $k''=1$ sur
$|\lambda^j,\mu^0|-|\lambda^0,\mu^0|$\\
si $k''=1$ vrai.\\
si $k''=k'$ sur $|\lambda^j,\mu^0|-|\lambda^0,\mu^0|$ alors $k''=k'=1$\\
alors $k''=1$ sur $|\lambda^j,\mu^0|-|\lambda^0,\mu^0|$ ce qui
donne Q=$(\lambda^i,\mu^n,k'')$
\begin{displaymath}
avec \qquad k'' = \left\{ \begin{array}{ll}
k & \textrm{ $|\lambda^0,\mu^0|$}\\
1 & \textrm{ $|\lambda^j,\mu^t|-|\lambda^0,\mu^0|$}
\end{array} \right.
\end{displaymath}
D'autre part,\\
$Q'=\big[(\mu^t,\mu^n,(a,b))(\lambda^j,\mu^t,k)\big](\lambda^i,\lambda^j,(c,d))$\\
=$(\lambda^j,\mu^n,k)(\lambda^i,\lambda^j,(c,d))$ \quad voir c)\\
=$(\lambda^i,\mu^n,k'')$  \quad voir d)\\
\begin{displaymath}
avec \qquad k' = \left\{ \begin{array}{ll}
k & \textrm{ $|\lambda^0,\mu^0|$}\\
1 & \textrm{ $|\lambda^0,\mu^t|-|\lambda^0,\mu^0|$}\\
k & \textrm{ $|\lambda^j,\mu^0|-|\lambda^0,\mu^0|$}\\
1 & \textrm{ $H'(\lambda^j,\mu^t)$}
\end{array} \right.
\end{displaymath}
et
\begin{displaymath}
avec \qquad k'' = \left\{ \begin{array}{ll}
k' & \textrm{ $|\lambda^0,\mu^0|$}\\
1 & \textrm{ $|\lambda^j,\mu^0|-|\lambda^0,\mu^0|$}\\
k' & \textrm{ $|\lambda^0,\mu^n|-|\lambda^0,\mu^0|$}\\
1 & \textrm{ $H'(\lambda^j,\mu^n)$}
\end{array} \right. \quad = \quad\left\{ \begin{array}{ll}
k & \textrm{ $|\lambda^0,\mu^0|$}\\
1 & \textrm{ $|\lambda^0,\mu^t|-|\lambda^0,\mu^0|$}\\
1 & \textrm{ $|\lambda^j,\mu^0|-|\lambda^0,\mu^0|$}\\
1 & \textrm{ $H'(\lambda^j,\mu^t)$}
\end{array} \right.
\end{displaymath}
On va verifier que k$''=1$ sur
$|\lambda^0,\mu^t|-|\lambda^0,\mu^0|$\\
On a k$''=k'=1$ ou $k''=1$ sur $|\lambda^0,\mu^t|-|\lambda^0,\mu^0|$\\
Donc $k''=1$ sur $|\lambda^0,\mu^t|-|\lambda^0,\mu^0|$.\\
Alors $Q'=(\lambda^i,\mu^n,k'')$
\begin{displaymath}
avec \qquad k'' = \left\{ \begin{array}{ll}
k & \textrm{ $|\lambda^0,\mu^0|$}\\
1 & \textrm{ $|\lambda^j,\mu^t|-|\lambda^0,\mu^0|$}
\end{array} \right.
\end{displaymath}
Donc $Q=Q'$ sur $|\lambda^j,\mu^t|$\\
alors
$\big[(\mu^t,\mu^n,(a,b))(\lambda^j,\mu^t,k)\big](\lambda^i,\lambda^j,(c,d))=(\mu^t,\mu^n,(a,b))\big[(\lambda^j,\mu^t,k)(\lambda^i,\lambda^j,(c,d))\big]$\\
\underline{\'equation 6}\\
\quad$\lambda^i\mapsto\lambda^j\mapsto\lambda^t\mapsto\mu^n$\\
Il y a deux cas:\\
1) $\lambda \in V$\\
2) $\lambda \in U$\\
Cas 1)\\
$Q=\big[(\lambda^0,\mu^n,k)(\lambda^0,\lambda^0,(a,b))\big](\lambda^i,\lambda^0,(c,d))$\\
=$(\lambda^0,\mu^n,1)(\lambda^i,\lambda^0,(c,d))$\quad car $\lambda \in V$\\
=$(\lambda^i,\mu^n,1)$\\
$Q'=(\lambda^0,\mu^n,k)\big[(\lambda^0,\lambda^0,(a,b))(\lambda^i,\lambda^0,(c,d))\big]$\\
=$(\lambda^0,\mu^n,k)(\lambda^i,\lambda^0,(c,b))$\\
=$(\lambda^i,\mu^n,1)$ \quad car $\lambda \in V$\\
Donc $Q=Q'$\\
Cas 2)\\
$$t=0,j=0$$
$Q=\big[(\lambda^0,\mu^n,k)(\lambda^0,\lambda^0,(a,b))\big](\lambda^i,\lambda^0,(c,d))$\\
=$(\lambda^0,\mu^n,k)(\lambda^i,\lambda^0,(c,d))$\quad voir d)\\
=$(\lambda^i,\mu^n,k)$ voir d)\\
D'autre part,\\
$Q'=(\lambda^0,\mu^n,k)\big[(\lambda^0,\lambda^0,(a,b))(\lambda^i,\lambda^0,(c,d))\big]$\\
=$(\lambda^0,\mu^n,k)(\lambda^i,\lambda^0,(c,b))$\\
=$(\lambda^i,\mu^n,k)$ \quad voir d)\\
Donc $Q=Q'$\\
$$t=0,j\neq 0$$
$Q=\big[(\lambda^0,\mu^n,k)(\lambda^j,\lambda^0,(a,b))\big](\lambda^i,\lambda^j,(c,d))$\\
=$(\lambda^j,\mu^n,k)(\lambda^i,\lambda^j,(c,d))$ voir d)\\
=$(\lambda^i,\mu^n,k')$ voir d)\\
\begin{displaymath}
avec \qquad k' = \left\{ \begin{array}{ll}
k & \textrm{ $|\lambda^0,\mu^0|$}\\
1 & \textrm{ $|\lambda^j,\mu^0|-|\lambda^0,\mu^0|$}\\
k & \textrm{ $|\lambda^0,\mu^n|-|\lambda^0,\mu^0|$}\\
1 & \textrm{ $H'(\lambda^j,\mu^n)$}\\
\end{array} \right.
\end{displaymath}
Donc Q=$(\lambda^i,\mu^n,k)$ sur $|\lambda^0,\mu^n|$\\
D$'autre$ part:\\
$Q'=(\lambda^0,\mu^n,k)\big[(\lambda^j,\lambda^0,(a,b))(\lambda^i,\lambda^j,(c,d))\big]$\\
=$(\lambda^0,\mu^n,k)(\lambda^i,\lambda^0,(c,b))$\\
=$(\lambda^i,\mu^n,k)$ voir d)\\
Alors $Q=Q'$ sur $|\lambda^0,\mu^n|$
$$t\neq 0,j=0$$
$Q=\big[(\lambda^t,\mu^n,k)(\lambda^0,\lambda^t,(a,b))\big](\lambda^i,\lambda^0,(c,d))$\\
=$(\lambda^0,\mu^n,k')(\lambda^i,\lambda^t,(c,b))$ voir d)\\
=$(\lambda^i,\mu^n,k')$ voir d)\\
\begin{displaymath}
avec \qquad k' = \left\{ \begin{array}{ll}
k & \textrm{ $|\lambda^0,\mu^0|$}\\
1 & \textrm{ $|\lambda^t,\mu^0|-|\lambda^0,\mu^0|$}\\
k & \textrm{ $|\lambda^0,\mu^n|-|\lambda^0,\mu^0|$}\\
1 & \textrm{ $H'(\lambda^t,\mu^n)$}\\
\end{array} \right.
\end{displaymath}
$Q'=(\lambda^t,\mu^n,k)\big[(\lambda^0,\lambda^t,(a,b))(\lambda^i,\lambda^0,(c,d))\big]$\\
=$(\lambda^t,\mu^n,k)(\lambda^i,\lambda^t,(c,b))$\\
=$(\lambda^i,\mu^n,k')$ voir d)\\
\begin{displaymath}
avec \qquad k' = \left\{ \begin{array}{ll}
k & \textrm{ $|\lambda^0,\mu^0|$}\\
1 & \textrm{ $|\lambda^t,\mu^0|-|\lambda^0,\mu^0|$}\\
k & \textrm{ $|\lambda^0,\mu^n|-|\lambda^0,\mu^0|$}\\
1 & \textrm{ $H'(\lambda^t,\mu^n)$}\\
\end{array} \right.
\end{displaymath}
Alors $Q=Q'$
$$t\neq 0, j\neq 0$$
$Q=\big[(\lambda^t,\mu^n,k)(\lambda^j,\lambda^t,(a,b))\big](\lambda^i,\lambda^j,(c,d))$\\
=$(\lambda^j,\mu^n,k')(\lambda^i,\lambda^j,(c,b))$ voir d)\\
=$(\lambda^i,\mu^n,k'')$ voir d)\\
\begin{displaymath}
avec \qquad k' = \left\{ \begin{array}{ll}
k & \textrm{ $|\lambda^0,\mu^0|$}\\
1 & \textrm{ $|\lambda^t,\mu^0|-|\lambda^0,\mu^0|$}\\
k & \textrm{ $|\lambda^0,\mu^n|-|\lambda^0,\mu^0|$}\\
1 & \textrm{ $H'(\lambda^t,\mu^n)$}\\
\end{array} \right.
\end{displaymath}
et
\begin{displaymath}
avec \qquad k'' = \left\{ \begin{array}{ll}
k' & \textrm{ $|\lambda^0,\mu^0|$}\\
1 & \textrm{ $|\lambda^j,\mu^0|-|\lambda^0,\mu^0|$}\\
k' & \textrm{ $|\lambda^0,\mu^n|-|\lambda^0,\mu^0|$}\\
1 & \textrm{ $H'(\lambda^j,\mu^n)$}\\
\end{array} \right. \quad = \quad\left\{ \begin{array}{ll}
k & \textrm{ $|\lambda^0,\mu^0|$}\\
1 & \textrm{ $|\lambda^t,\mu^0|-|\lambda^0,\mu^0|$}\\
k & \textrm{ $|\lambda^0,\mu^n|-|\lambda^0,\mu^0|$}\\
1 & \textrm{ $H'(\lambda^t,\mu^n)$}\\
\end{array} \right.
\end{displaymath}
D'autre part,\\
$Q'=(\lambda^t,\mu^n,k)\big[(\lambda^j,\lambda^t,(a,b))(\lambda^i,\lambda^j,(c,d))\big]$\\
=$(\lambda^t,\mu^n,k)(\lambda^i,\lambda^t,(c,b))$\\
=$(\lambda^i,\mu^n,k')$ voir d)\\
\begin{displaymath}
avec \qquad k' = \left\{ \begin{array}{ll}
k & \textrm{ $|\lambda^0,\mu^0|$}\\
1 & \textrm{ $|\lambda^t,\mu^0|-|\lambda^0,\mu^0|$}\\
k & \textrm{ $|\lambda^0,\mu^n|-|\lambda^0,\mu^0|$}\\
1 & \textrm{ $H'(\lambda^t,\mu^n)$}\\
\end{array} \right.
\end{displaymath}
Alors $Q=Q'$ sur $|\lambda^t,\mu^n|$\\
Donc l'\'equation 6 est associative.\\
\underline{\'equation 7} est comme \underline{\'equation 6}\\
\underline{\'equation 8}\\
\quad$\lambda^i\mapsto\lambda^j\mapsto\lambda^t\mapsto\lambda^n$\\
Q=$\big[(\lambda^t,\lambda^n,(a,b))(\lambda^j,\lambda^t,(c,d))\big](\lambda^i,\lambda^j,(e,f))$=\\
$(\lambda^j,\lambda^n,(c,b))(\lambda^i,\lambda^j,(e,f))$=\\
$(\lambda^i,\lambda^n,(e,b))$\\
Q$'=(\lambda^t,\lambda^n,(a,b))\big[(\lambda^j,\lambda^t,(c,d))(\lambda^i,\lambda^j,(e,f))\big]=$\\
$(\lambda^t,\lambda^n,(a,b))(\lambda^i,\lambda^t,(e,d))=$\\
$(\lambda^i,\lambda^n,(e,b))$\\
Donc $Q=Q'$\\
On a verifier toutes les \'equations donc A est bien une
cat\'egorie
associ\'ee \`a M donc $Cat(M)\neq \emptyset$.\\

\begin{corollary}:Soit M=$(m_{ij})$ une matrice carr\'ee d'entiers, alors
 $Cat(M)\neq \emptyset$ si et seulement si, pour toute sous-matrice $M'$ de $M$, de taille $\leq 4$, $Cat (M')\neq\emptyset$.
\end{corollary}
Preuve:\\
$\Rightarrow)$ evidente.\\
$\Leftarrow)$ La condition d'acceptabilit\'e, et les in\'egalit\'es du th\'eor\`eme
\ref{gen1}, font intervenir $\leq 4$ objets \`a la fois.


\begin{thebibliography}{MM}


\bibitem{Allouch}
S. Allouch. Classification des cat\'egories
finies. http://math.unice.fr/~carlos/documents/allouchJun07.pdf,
M\'emoire de M2, 15 juin(2007).

\bibitem{Allouch2}
S. Allouch. Sur l'existence d'une cat\'egorie ayant une matrice 
strictement positive donn\'ee. Preprint arXiv:0806.0060v1 (2008). 

\bibitem{BergerLeinster}
C. Berger, T. Leinster. The Euler characteristic of a category as
the sum of a divergent series. Homology, Homotopy Appl. 10 (2008),
41-51.

\bibitem{BrinkmannMcKay}
G. Brinkmann, B. McKay. Posets on up to 16 Points. Order 19(2002),
147-179.

\bibitem{CuntzHeckenberger}
M. Cuntz, I. Heckenberger. Weyl groupoids with at most three
objects. Preprint arXiv:0805.1810v1 [math.GR].

\bibitem{FioreLuckSauer}
T. Fiore, W. L\"uck, R. Sauer. 
Finiteness obstructions and Euler characteristics of categories.
Preprint arXiv:0908.3417 (2009). 


\bibitem{FlemingGuntherRosebrugh}
M. Fleming, R. Gunther, R. Rosebrugh. A Database of Categories.
Journal of Symbolic Computation 35 (2003), 127-135.

\bibitem{ForresterBarker}
M. Forrester-Barker. Group objects and internal categories.
Preprint arXiv:math/0212065v1.

\bibitem{Kapranov}
M. Kapranov. On the derived categories of coherent sheaves on some
homogeneous spaces. Invent. Math. 92 (1988), 479-508.


\bibitem{Leinster}
T. Leinster. The Euler characteristic of a category. Doc. Math., 13 (2008), 21-49,
arXiv:math/0610260v1.


\end{thebibliography}
\end{document}